\documentclass[12pt,twoside]{amsart}
\usepackage{amssymb,amsmath,amsthm, amscd, enumerate, mathrsfs}
\usepackage{graphicx, hhline}
\usepackage[all]{xy}
\usepackage{hyperref}
\usepackage{geometry}
\title[On mixed-$\omega$-sheaves] 
{On mixed-$\omega$-sheaves} 
\author{Osamu Fujino} 
\date{2023/9/8, version 0.66}
\address{Department of 
Mathematics, Graduate School of Science, 
Kyoto University, Kyoto 606-8502, Japan}
\email{fujino@math.kyoto-u.ac.jp}

\newcommand{\Spec}[0]{{\operatorname{Spec}}}
\newcommand{\rank}[0]{{\operatorname{rank}}}
\newcommand{\Exc}[0]{{\operatorname{Exc}}}
\newcommand{\Supp}[0]{{\operatorname{Supp}}}
\newtheorem{thm}{Theorem}[section]
\newtheorem{lem}[thm]{Lemma}
\newtheorem{cor}[thm]{Corollary}

\newtheorem{conj}[thm]{Conjecture}
\newtheorem*{claim}{Claim}

\theoremstyle{definition}
\newtheorem{defn}[thm]{Definition}
\newtheorem{rem}[thm]{Remark}
\newtheorem{ex}[thm]{Example}
\newtheorem{step}{Step}
\newtheorem*{ack}{Acknowledgments} 

\makeatletter
    
    \@addtoreset{equation}{section}
\makeatother

\begin{document}

\maketitle 

\begin{abstract}
We introduce the notion of mixed-$\omega$-sheaves and 
use it for the study of a relative version of Fujita's freeness conjecture. 
It is related to the Iitaka conjecture. 
We note that the notion of mixed-$\omega$-sheaves 
is a generalization of that of Nakayama's $\omega$-sheaves in some sense. 
One of the main motivations of this paper is to make 
Nakayama's theory of $\omega$-sheaves more 
accessible and make it applicable to the study of log canonical pairs. 
\end{abstract}

\tableofcontents 

\section{Introduction}\label{q-sec1} 

Let us recall Fujita's famous freeness conjecture on adjoint bundles. 
Note that everything is defined over $\mathbb C$, the complex number 
field, in this paper. 

\begin{conj}[Takao Fujita]\label{q-conj1.1}
Let $X$ be a smooth projective 
variety with $\dim X=n$ and let 
$\mathcal L$ be any ample invertible sheaf on $X$. 
Then $\omega_X\otimes \mathcal L^{\otimes l}$ is 
generated by global sections for every $l\geq n+1$. 
\end{conj}

Although there have already been many related 
results, Conjecture \ref{q-conj1.1} 
is still open. 
As a generalization of Conjecture \ref{q-conj1.1}, 
Popa and Schnell proposed the following conjecture, 
which is a relative version of Fujita's freeness conjecture. 

\begin{conj}[{Popa and Schnell, see 
\cite[Conjecture 1.3]{popa-schnell}}]\label{q-conj1.2}
Let $f\colon X\to Y$ be a surjective morphism between smooth 
projective varieties with $\dim Y=n$. 
Let $\mathcal L$ be any ample invertible sheaf on $Y$. 
Then, for every positive integer $k$, the sheaf 
$$
f_*\omega^{\otimes k}_X\otimes \mathcal L^{\otimes l} 
$$ 
is generated by global sections for $l\geq k(n+1)$. 
\end{conj}

We can find some interesting results on Conjecture 
\ref{q-conj1.2} in \cite{denf}, \cite{dutta}, \cite{dutta-murayama}, 
and \cite{iwai}. 
In this paper, we do not directly treat Conjecture \ref{q-conj1.2}. 
When the base space $Y$ is a curve in Conjecture 
\ref{q-conj1.2}, 
we have the following stronger result. 
We note that we are mainly interested in 
the case where $k\geq 2$ for some geometric applications. 

\begin{thm}\label{q-thm1.3}
Let $f\colon X\to C$ be a surjective morphism 
from a smooth projective variety $X$ 
onto a smooth projective curve $C$. 
Let $\mathcal H$ be an ample invertible sheaf on $C$ 
with $\deg \mathcal H\geq 2$ and 
let $k$ be any positive integer. 
Then the sheaf $$
f_*\mathcal \omega^{\otimes k}_{X/C}
\otimes \omega_C\otimes \mathcal H
$$ 
is generated by global sections. 
Let $\mathcal L$ be an ample invertible sheaf on $C$. 
Then 
$$
f_*\omega^{\otimes k}_X\otimes \mathcal L^{\otimes l} 
$$ 
is generated by global sections for $l\geq 2k$. 
In particular, Conjecture \ref{q-conj1.2} holds true when 
the base space is a smooth projective curve. 
\end{thm}

Here, we give a detailed proof of Theorem \ref{q-thm1.3} in order 
to explain our idea. 

\begin{proof}[Proof of Theorem \ref{q-thm1.3}]
If $f_*\omega^{\otimes k}_{X/C}=0$, 
then there is nothing to prove. 
So we assume that $f_*\omega^{\otimes k}_{X/C}
\ne 0$. 
We take any closed point $P$. 
\begin{claim}\label{q-claim}
$H^1(C, f_*\omega^{\otimes k}_{X/C}\otimes \omega_C
\otimes \mathcal H\otimes \mathcal O_C(-P))=0$. 
\end{claim} 
\begin{proof}[Proof of Claim] 
By \cite[Theorem 1]{kawamata-curve}, 
$f_*\omega^{\otimes k}_{X/C}$ is a nef locally free sheaf. 
This fact also follows from 
Viehweg's weak positivity theorem since $C$ is a smooth 
projective curve. 
Therefore, $\mathcal E:=f_*\omega^{\otimes k}_{X/C}\otimes 
\mathcal H\otimes \mathcal O_C(-P)$ is ample. 
If $H^1(C, \mathcal E\otimes \omega_C)\ne 0$, 
then we get $H^0(C, \mathcal E^*)\ne 0$ by Serre duality. 
This implies that 
there is a nontrivial inclusion 
$0\to \mathcal O_C\to\mathcal E^*$. 
By taking the dual of this inclusion, we have the following 
surjection $\mathcal E\to \mathcal O_C\to 0$. 
This is a contradiction since $\mathcal E$ is ample. 
Hence we have $H^1(C, \mathcal E\otimes \omega_C)=0$. 
\end{proof}
By the Claim, the natural restriction map 
$$
H^0(C, f_*\omega^{\otimes k}_{X/C}\otimes 
\omega_C\otimes \mathcal H)\to 
f_*\omega^{\otimes k}_{X/C}\otimes 
\omega_C\otimes \mathcal H\otimes \mathbb C(P)
$$
is surjective. 
This means that 
$f_*\omega^{\otimes k}_{X/C}\otimes 
\omega_C\otimes \mathcal H$ is generated by global sections. 
Since $C$ is a smooth projective curve, 
$\omega_C\otimes \mathcal L^{\otimes 2}$ 
is generated by global sections. 
This implies that 
$$
f_*\omega^{\otimes k}_X\otimes \mathcal L^{\otimes l} 
\simeq f_*\omega^{\otimes k}_{X/C}\otimes 
\omega_C\otimes \mathcal L^{\otimes (l-2(k-1))}\otimes 
\left(\omega_C\otimes \mathcal L^{\otimes 2}\right)
^{\otimes (k-1)}
$$ 
is generated by global sections because 
$\deg \mathcal L^{\otimes (l-2(k-1))}\geq 2$ by assumption. 
\end{proof}

A key point of Theorem \ref{q-thm1.3} is 
the fact that $f_*\omega^{\otimes k}_{X/C}$ is a nef locally 
free 
sheaf on $C$ for 
every positive integer $k$. 
When $f\colon X\to Y$ is a weakly semistable morphism 
in the sense of Abramovich--Karu (see \cite{abramovich-karu}), it is 
conjectured that $f_*\omega^{\otimes k}_{X/Y}$ is a nef locally 
free sheaf on $Y$ for every positive integer $k$ 
(see \cite{fujino-direct} and 
\cite[Conjecture 3.14]{fujino-zucker}). 
Hence it is natural to consider the case where 
$f\colon X\to Y$ is weakly semistable. 

\begin{thm}[see 
Theorem \ref{p-thm8.2}]\label{q-thm1.4}
Let $f\colon X\to Y$ be a surjective 
morphism from a normal projective variety $X$ onto a smooth projective 
variety $Y$ with connected fibers. 
Assume that $f$ is weakly semistable in the sense of 
Abramovich--Karu and that the geometric generic 
fiber $X_{\overline \eta}$ of $f\colon  X\to Y$ has a good minimal model. 
Let $H$ be an ample Cartier divisor on $Y$,  
let $k$ be a positive integer with $k\geq 2$, and let 
$A$ be an ample Cartier divisor on $Y$ such that 
$|A|$ is free. Then 
$$
\left(\bigotimes ^s f_*\omega^{\otimes k}_{X/Y}\right)\otimes \omega_Y 
\otimes \mathcal O_Y(H+nA)
$$ 
is generated by global sections for all integers $s\geq 1$, 
where $n=\dim Y$. 
\end{thm}

By Theorems \ref{q-thm1.3} and 
\ref{q-thm1.4}, we propose a new 
conjecture similar to Conjecture \ref{q-conj1.2}. 

\begin{conj}\label{q-conj1.5}
Let $f\colon  X\to Y$ be a surjective morphism 
between smooth projective varieties with $\dim Y=n$. 
Let $\mathcal L$ be any ample invertible sheaf on $Y$. 
Then, for every positive integer $k$, the sheaf 
$$
f_*\omega^{\otimes k}_{X/Y}\otimes 
\omega_Y\otimes \mathcal L^{\otimes l}
$$ 
is generically generated by global sections for $l\geq n+1$. 
More precisely, 
$$
f_*\omega^{\otimes k}_{X/Y}\otimes 
\omega_Y\otimes \mathcal L^{\otimes l}
$$ 
is generated by global section on $U$ for $l\geq n+1$, where 
$U$ is the largest Zariski open set of $Y$ such that 
$f$ is smooth over $U$.  
\end{conj}

Even when $f$ is the identity in Conjecture \ref{q-conj1.5}, 
the generically generation, 
which is a special case of Fujita's conjecture (see Conjecture 
\ref{q-conj1.1}), looks like a difficult open problem. 
We note that an example (see Example \ref{p-ex10.1}) shows that 
$f_*\omega^{\otimes k}_{X/Y}\otimes 
\omega_Y \otimes \mathcal L^{\otimes l}$ is not always 
generated by global sections in Conjecture \ref{q-conj1.5}. 
We only expect that it is generically generated by global sections. 
The best known result on Conjecture \ref{q-conj1.5} 
is the following theorem. 
The author learned it from Masataka Iwai. 

\begin{thm}[Masataka Iwai, see Theorem \ref{p-thm12.1}]\label{q-thm1.6} 
Let $f\colon X\to Y$ be a surjective morphism between 
smooth projective varieties with connected fibers
and let $\mathcal L$ be an ample invertible sheaf on $Y$. 
Let $U$ be the largest Zariski open set of $Y$ such that 
$f$ is smooth over $U$. 
We put $\dim Y=n$. 
Then 
$$
f_*\omega^{\otimes a}_{X/Y}\otimes \omega_Y\otimes \mathcal L^{\otimes b}
$$ 
is generated by global sections on $U$ for all integers 
$a\geq 1$ and $b\geq \frac{n(n+1)}{2}+1$. 
\end{thm}

The proof of Theorem \ref{q-thm1.6} 
is essentially analytic. 
We will give a sketch of the proof of Theorem \ref{q-thm1.6} 
in Section \ref{p-sec12}. On the other hand, 
by Nakayama's theory of $\omega$-sheaves 
(see \cite[Chapter V]{nakayama}), 
we can prove: 

\begin{thm}\label{q-thm1.7}
Let $f\colon X\to Y$ be a surjective morphism between smooth 
projective varieties and let $H$ be an ample divisor 
on $Y$ such that 
$|H|$ is free. 
We put $\dim Y=n$. 
Then 
$$
\left(\bigotimes ^s f_*\omega^{\otimes k}_{X/Y}\right)^{**} 
\otimes \omega_Y\otimes \mathcal O_Y(lH)
$$ 
is generically generated by global sections for all integers $k\geq 1$, 
$s\geq 1$, and $l\geq n+1$. 

Let $H^\dag$ be an ample divisor on $Y$ 
such that $|H^\dag|$ is not necessarily free. 
Then 
the sheaf 
$$
\left(\bigotimes ^s f_*\omega^{\otimes k}_{X/Y}\right)^{**} 
\otimes \omega_Y\otimes \mathcal O_Y(lH^\dag)
$$ 
is generically generated by global sections for all integers $k\geq 1$, 
$s\geq 1$, and $l\geq n^2+\min\{2, k\}$. 
\end{thm}

In this paper, we will introduce the notion of mixed-$\omega$-sheaves, 
which is a generalization of that of Nakayama's $\omega$-sheaves,  
and establish 
Theorems \ref{q-thm1.8} and \ref{q-thm1.9}. 
Then we will obtain  
Theorem \ref{q-thm1.7} as a special case of 
Theorems \ref{q-thm1.8} and \ref{q-thm1.9}. 

\begin{thm}\label{q-thm1.8}
Let $f\colon X\to Y$ be a surjective morphism 
from a normal projective variety $X$ onto a smooth 
projective variety $Y$ with $\dim Y=n$. 
Let $\Delta$ be an effective $\mathbb R$-divisor 
on $X$ such that $K_X+\Delta$ is $\mathbb R$-Cartier 
and that $(X, \Delta)$ is 
log canonical over a nonempty Zariski open set of $Y$. 
Let $L$ be a Cartier divisor on $X$ with 
$L\sim _{\mathbb R}k(K_{X/Y}+\Delta)$ for some 
positive integer $k$.  
Let $H$ be a big Cartier divisor on $Y$ such that 
$|H|$ is free. 
Then 
$$
\left(\bigotimes ^s f_*\mathcal O_X(L)\right)^{**} 
\otimes \mathcal O_Y(K_Y+lH)
$$ 
is generically generated by global sections for all integers $s\geq 1$ 
and $l\geq n+1$. 
\end{thm}

\begin{thm}\label{q-thm1.9}
In Theorem \ref{q-thm1.8}, 
we assume that $H^\dag$ is a nef and big 
Cartier divisor 
on $Y$ such that $|H^\dag|$ is not necessarily 
free. 
Then we have that 
$$ 
\left(\bigotimes ^s f_*\mathcal O_X(L)\right)^{**} 
\otimes \mathcal O_Y(K_Y+lH^\dag)
$$ 
is generically generated by global sections 
for all integers $s\geq 1$ and $l \geq n^2+\min\{2, k\}$. 
\end{thm} 

We note that Iwai's analytic method can not be applied to 
log canonical pairs 
because it depends on $L^2$ method. 
We also note that Nakayama's theory 
of $\omega$-sheaves can not be directly applied 
to the study of log canonical pairs. 

\medskip 

Let us quickly explain the idea of the 
proof of Theorem \ref{q-thm1.7}, which is 
mainly due to Nakayama (see \cite[Chapter V]{nakayama}). 
Let $f\colon X\to Y$ be a surjective morphism between 
smooth projective 
varieties and let $H$ be an ample Cartier divisor on $Y$ such that 
$|H|$ is free. 
We fix a positive integer $k\geq 2$. 
Then we can construct a surjective morphism 
$g\colon Z\to Y$ from a smooth projective 
variety $Z$ and a direct summand $\mathcal F$ 
of $g_*\mathcal O_Z(K_Z)$ such that 
there exists a generically isomorphic 
injection 
$$
\mathcal F\hookrightarrow 
\left(f_*\omega^{\otimes k}_{X/Y}\otimes 
\omega_Y\otimes \mathcal O_Y(H)\right)^{**}. 
$$ 
By Koll\'ar's vanishing theorem, 
we have 
$$
H^i(Y, \mathcal F\otimes \mathcal O_Y((n+1-i)H))=0
$$ 
for every $i>0$, where 
$n=\dim Y$. 
Therefore, by Castelnuovo--Mumford regularity, 
$\mathcal F\otimes \mathcal O_Y((n+1)H)$ is generated 
by global sections. 
This implies that 
$$
\left(f_*\omega^{\otimes k}_{X/Y}\right)^{**}\otimes 
\omega_Y\otimes \mathcal O_Y((n+2)H)
$$
is generically generated by global sections. 
Note that we do not try to 
establish any vanishing theorems for 
$f_*\omega^{\otimes k}_{X/Y}\otimes \omega_Y\otimes 
\mathcal O_Y(H)$ directly. 
By the above observation, it is natural to consider: 

\begin{defn}[Mixed-$\omega$-sheaf and 
pure-$\omega$-sheaf, see Definition \ref{p-def5.1}]\label{q-def1.10} 
A torsion-free coherent sheaf $\mathcal F$ on a normal quasi-projective 
variety $W$ is called a {\em{mixed-$\omega$-sheaf}} 
if there exist a projective surjective morphism 
$f\colon V\to W$ from a smooth quasi-projective variety $V$ and a simple 
normal crossing divisor $D$ on $V$ such that 
$\mathcal F$ is a direct summand of $f_*\mathcal O_V(K_V+D)$. 
When $D=0$, $\mathcal F$ is called a {\em{pure-$\omega$-sheaf}} on $W$. 
\end{defn}

For the study of klt pairs, 
the notion of pure-$\omega$-sheaves is sufficient and is essentially 
due to Nakayama (see \cite[Chapter V]{nakayama}). 
In this paper, we study some basic properties 
of mixed-$\omega$-sheaves. 
They are indispensable for the study of 
log canonical pairs. 
Of course, the theory of 
mixed-$\omega$-sheaves (resp.~pure-$\omega$-sheaves) 
in this paper is based on that of mixed (resp.~pure) Hodge 
structures. 
Roughly speaking, Nakayama only treats pure-$\omega$-sheaves 
in \cite[Chapter V]{nakayama}. However, his theory of 
$\omega$-sheaves is more sophisticated and 
some of his results are much sharper than 
ours. We do not try to make the framework 
discussed in this paper supersede Nakayama's theory 
of $\omega$-sheaves in \cite[Chapter V]{nakayama}. 
One of the main purposes of this 
paper is to make Nakayama's theory of $\omega$-sheaves 
more accessible and 
make it applicable to the study of log canonical pairs. 
Theorem \ref{p-thm9.3} (and Remark 
\ref{p-rem9.4}) is one of the main results of this 
paper, which we call a fundamental theorem of the theory 
of mixed-$\omega$-sheaves. 

\begin{thm}[{see \cite[Chapter V, 3.35.~Theorem]{nakayama}, 
Theorem \ref{p-thm9.3}, and Remark \ref{p-rem9.4}}]\label{q-thm1.11}
Let $f\colon X\to Y$ be a surjective morphism 
from a normal projective variety $X$ onto a smooth projective 
variety $Y$. 
Let $L$ be a Cartier divisor 
on $X$ and let $\Delta$ be an effective $\mathbb R$-divisor 
on $X$ such that 
$K_X+\Delta$ is $\mathbb R$-Cartier. 
Let $D$ be an $\mathbb R$-divisor on $Y$. 
Let $k$ be a positive integer with $k\geq 2$. 
Assume the following conditions: 
\begin{itemize}
\item[(i)] $(X, \Delta)$ is log canonical {\em{(}}resp.~klt{\em{)}} 
over a nonempty Zariski 
open set of $Y$, 
and 
\item[(ii)] $L+f^*D-k(K_{X/Y}+\Delta)-f^*A$ is 
semi-ample for some big $\mathbb R$-divisor 
$A$ on $Y$. 
\end{itemize}
If $f_*\mathcal O_Y(L)\ne 0$, then there exist a mixed-$\omega$-big-sheaf 
{\em{(}}resp.~pure-$\omega$-big-sheaf{\em{)}} 
$\mathcal F$ on $Y$ and a generically 
isomorphic injection 
$$
\mathcal F\hookrightarrow 
\mathcal O_Y(K_Y+\lceil D\rceil)\otimes \left(f_*\mathcal O_X(L)\right)^{**}. 
$$ 
\end{thm}

For the precise definition of mixed-$\omega$-big-sheaves and 
pure-$\omega$-big-sheaves, see Definition \ref{p-def5.3} below. 

As an application of Theorem \ref{q-thm1.11}, 
we give a detailed proof of: 

\begin{thm}[{\cite[Chapter V, 4.1.~Theorem (1)]{nakayama},  
\cite[Section 3]{fujino-corri}, and Theorem \ref{p-thm11.3}}]\label{q-thm1.12} 
Let $f\colon X\to Y$ be a surjective morphism from a normal 
projective variety $X$ onto a smooth projective variety $Y$ 
with connected fibers. 
Let $\Delta$ be an effective $\mathbb R$-divisor 
on $X$ such that 
$K_X+\Delta$ is $\mathbb R$-Cartier and that 
$(X, \Delta)$ is log canonical over a nonempty Zariski open set 
of $Y$. 
Let $D$ be an $\mathbb R$-Cartier $\mathbb R$-divisor 
on $X$ such that $D-(K_{X/Y}+\Delta)$ is nef. 
Then, for any $\mathbb R$-divisor $Q$ on $Y$, 
we have 
$$
\kappa _\sigma(X, D+f^*Q)\geq \kappa _\sigma(F, D|_F)
+\kappa (Y, Q)
$$ 
and 
$$
\kappa _\sigma(X, D+f^*Q)\geq \kappa (F, D|_F)
+\kappa _\sigma(Y, Q)
$$ 
where $F$ is a sufficiently general fiber of $f\colon X\to Y$. 
\end{thm}

We note that $\kappa_\sigma(X, D)$ and $\kappa (X, D)$ 
denote Nakayama's numerical dimension and 
the Iitaka dimension of $D$, respectively. 
Theorem \ref{q-thm1.12} already played a crucial role in the theory 
of minimal models. 
We need Theorem \ref{q-thm1.12} for the proof of the following 
famous and fundamental 
result on the existence theorem of good minimal models 
for klt pairs. 

\begin{thm}[{\cite[Remark 2.6]{dhp}, 
\cite[Theorem 4.3]{gongyo-lehmann}, and 
\cite[Theorem 3.2]{fujino-corri}}]\label{q-thm1.13} 
Let $(X, \Delta)$ be a projective 
klt pair such that $\Delta$ is a $\mathbb Q$-divisor. 
Then $(X, \Delta)$ has a good 
minimal model if and only if 
$\kappa_\sigma(X, K_X+\Delta)=
\kappa (X, K_X+\Delta)$. 
\end{thm}

We explain the organization of this paper. 
In Section \ref{p-sec2}, we collect some basic definitions. 
In Section \ref{p-sec3}, we prepare some useful and important lemmas. 
They will play a crucial role in this paper. 
In Section \ref{p-sec4}, we quickly explain some basic properties 
of Viehweg's weakly positive sheaves and big sheaves. 
In Section \ref{p-sec5}, we introduce mixed-$\omega$-sheaves and 
mixed-$\omega$-big-sheaves. 
In Sections \ref{p-sec6} and \ref{p-sec7}, we prove 
some basic properties of mixed-$\omega$-sheaves based on 
the theory of mixed Hodge structures. 
In Section \ref{p-sec8}, we treat a very special but interesting case. 
More precisely, we treat weakly semistable morphisms 
$f\colon X\to Y$ in the sense of Abramovich--Karu with 
the assumption that the geometric generic fiber $X_{\overline \eta}$ 
of $f\colon X\to Y$ has 
a good minimal model. 
In this case, we can prove some strong results with 
the aid of the theory of minimal models. 
Section \ref{p-sec9} is the main part of this paper. 
We prove Theorem \ref{q-thm1.11}. 
Section \ref{p-sec10} is devoted to the proof of 
Theorems \ref{q-thm1.7}, \ref{q-thm1.8}, and \ref{q-thm1.9}. 
In Section \ref{p-sec11}, we prove Theorem \ref{q-thm1.12}, 
which has already played a crucial role 
in the theory of minimal models, and 
a slight generalization of the twisted weak positivity theorem, 
which is an extension of results of various people. 
In the final section:~Section \ref{p-sec12}, 
we give a sketch of the proof of Theorem \ref{q-thm1.6}. 
The proof is essentially analytic and is independent of 
the other sections. 

\begin{ack} 
During the preparation of this paper,
the author was partially supported by 
JSPS KAKENHI Grant Numbers JP16H03925, JP16H06337, 
JP19H01787, JP20H00111, 
JP21H00974, JP21H04994. 
He thanks Masataka Iwai, Takumi Murayama, and 
Sho Ejiri 
for useful comments, suggestions, and discussions. 
In 2020, he gave a series of lectures based on 
a preprint version of this paper in Tianyuan 
Workshop, which was organized by Lei Zhang. 
The author thanks him and every audience. 
\end{ack}

We will work over $\mathbb C$, the complex number field, throughout 
this paper. We note that a {\em{scheme}} is a separated 
scheme of finite type over $\mathbb C$ and 
a {\em{variety}} is an integral scheme. 

\section{Preliminary}\label{p-sec2}

In this section, we collect some basic definitions. 
For the details, see \cite{fujino-fundamental}, \cite{fujino-foundations}, and 
\cite{fujino-iitaka}. 

\medskip 
 
Let us start with the definition of canonical sheaves and canonical divisors. 

\begin{defn}[Canonical sheaf and canonical divisor]\label{p-def2.1}
Let $X$ be an equidimensional 
scheme of dimension $n$ and 
let $a\colon X\to \Spec \mathbb C$ be 
the structure morphism. 
Then the 
dualizing complex of $X$ is $\omega^\bullet_X
=a^!\mathbb C$, where $a^!$ is the functor 
obtained in \cite[Chapter VII.~Corollary 3.4 (a)]{hartshorne} 
(see also \cite[Section 3.3]{conrad}). 
Then we put 
$$
\omega_X:=h^{-n}(\omega^\bullet_X)
$$ 
and call it the {\em{canonical sheaf}} of $X$. 

We further assume that $X$ is normal. 
Then a {\em{canonical divisor}} $K_X$ of $X$ 
is a Weil divisor on $X$ such that 
$$
\mathcal O_{X_{\mathrm{sm}}}(K_X)\simeq 
\Omega^n_{X_{\mathrm{sm}}}
$$
holds, where $X_{\mathrm{sm}}$ is the largest smooth 
Zariski open set 
of $X$. 

It is well known that 
$$
\mathcal O_X(K_X)\simeq \omega_X
$$ 
holds when $X$ is normal. 

If $f\colon X\to Y$ is a morphism between Gorenstein schemes, 
then we put 
$$
\omega_{X/Y}:=\omega_X\otimes f^*\omega^{\otimes -1}_Y. 
$$

If $f\colon X\to Y$ is a morphism from a normal scheme $X$ to 
a normal Gorenstein scheme $Y$, 
then we put 
$$
K_{X/Y}:=K_X-f^*K_Y. 
$$
\end{defn}

Let us quickly look at 
the definition of singularities of pairs. 

\begin{defn}[Singularities of pairs]\label{p-def2.2}
Let $X$ be a normal variety and let $\Delta$ be 
an effective $\mathbb R$-divisor on $X$ such that 
$K_X+\Delta$ is $\mathbb R$-Cartier. 
Let $f\colon Y\to X$ be a resolution of singularities of $X$ such that 
$\Exc(f)\cup f^{-1}_*\Delta$ has simple normal crossing support, 
where $\Exc(f)$ is the exceptional locus of 
$f$ on $Y$ and $f^{-1}_*\Delta$ is the strict transform of 
$\Delta$ on $Y$. 
Then we can write 
$$
K_Y=f^*(K_X+\Delta)+\sum _i a_i E_i
$$ 
with $f_*\left(\sum _i a_i E_i\right)=-\Delta$. 
We say that $(X, \Delta)$ is {\em{log canonical}} 
(resp.~{\em{klt}}) if $a_i\geq -1$ (resp.~$a_i>-1$) for 
every $i$. If $\Delta=0$ and $a_i\geq 0$ holds for every $i$, 
then we say that $X$ has only {\em{canonical singularities}}. 

If $(X, \Delta)$ is log canonical and there exist a resolution 
of singularities $f\colon Y\to X$ as above and a prime 
divisor $E_i$ on $Y$ with $a_i=-1$, then 
$f(E_i)$ is called a {\em{log canonical center}} of 
$(X, \Delta)$. 
\end{defn}

\begin{defn}[Dlt pairs]\label{p-def2.3}
Let $(X, \Delta)$ be a log canonical pair. 
If there exists a resolution of singularities 
$f\colon Y\to X$ such that 
the exceptional locus $\Exc(f)$ of $f$ 
is a divisor on $Y$, 
$\Exc(f)\cup f^{-1}_*\Delta$ has simple normal 
crossing support, and 
$$
K_Y=f^*(K_X+\Delta)+\sum _i a_i E_i
$$
with $a_i>-1$ for every $f$-exceptional divisor $E_i$, 
then the pair $(X, \Delta)$ is called a {\em{dlt}} 
pair. 
\end{defn}

The following definitions are very useful 
in this paper. 

\begin{defn}[Horizontal and vertical 
divisors]\label{p-def2.4} 
Let $f\colon X\to Y$ be a dominant morphism between normal  
varieties and let $D$ be an $\mathbb R$-divisor on $X$. 
We can write $$D=D_{\mathrm{hor}}+D_{\mathrm{ver}}$$ such 
that every irreducible component of $D_{\mathrm{hor}}$ (resp.~
$D_{\mathrm{ver}}$) is mapped (resp.~not mapped) onto 
$Y$. 
If $D=D_{\mathrm{hor}}$ (resp.~$D=D_{\mathrm{ver}}$), then 
$D$ 
is said to be {\em{horizontal}} (resp.~{\em{vertical}}). 
\end{defn}

\begin{defn}[Operations for $\mathbb R$-divisors]\label{p-def2.5}
Let $D=\sum _i d_i D_i$ be an $\mathbb R$-divisor on a normal 
variety $X$, where 
$D_i$ is a prime divisor on $X$ for every $i$, 
$D_i\ne D_j$ for $i\ne j$, and $d_i \in 
\mathbb R$ for every $i$.   
Then we put $$\lfloor D\rfloor =\sum _i \lfloor d_i \rfloor D_i, 
\quad \{D\}=D-\lfloor 
D\rfloor, \quad \text{and} \quad\lceil D\rceil =-\lfloor -D\rfloor. $$ 
Note that $\lfloor d_i \rfloor$ is the integer which satisfies 
$d_i-1<\lfloor d_i \rfloor \leq d_i$. 
We also note that 
$\lfloor D\rfloor$, $\lceil D\rceil$, and $\{D\}$ are called 
the {\em{round-down}}, {\em{round-up}}, 
and {\em{fractional part}} of $D$ 
respectively. 

If $0\leq d_i \leq 1$ for every $i$, then 
we say that $D$ is a {\em{boundary}} 
$\mathbb R$-divisor on $X$. 
We note that $\sim _{\mathbb Q}$ (resp.~$\sim_{\mathbb R}$) 
denotes the {\em{$\mathbb Q$-linear}} 
(resp.~{\em{$\mathbb R$-linear}}) {\em{equivalence}} 
of $\mathbb Q$-divisors (resp.~$\mathbb R$-divisors). 

In this paper, we will repeatedly use 
the following notation:  
$$
D^{=1}=\sum _{d_i=1} D_i, \quad 
D^{>1}=\sum _{d_i>1} d_i D_i, \quad \text{and} \quad 
D^{<0}=\sum _{d_i<0}d_i D_i. 
$$
\end{defn}

We recall the following definition for the 
reader's convenience. 

\begin{defn}[Generic generation]\label{p-def2.6} 
Let $\mathcal F$ be a coherent 
sheaf on a variety $X$. 
We say that $\mathcal F$ is {\em{generated by global sections 
on}} $U$, where $U$ is a Zariski open set of $X$, 
if the natural map 
$$
H^0(X, \mathcal F)\otimes \mathcal O_X\to \mathcal F
$$ 
is surjective on $U$. 
We simply say that $\mathcal F$ is {\em{generated 
by global sections}} when $U=X$. 
We say that $\mathcal F$ is {\em{generically 
generated by global sections}} if $\mathcal F$ is 
generated by global sections on some nonempty 
Zariski open set of $X$. 
\end{defn}

We close this section with the definition of 
exceptional divisors for proper 
surjective morphisms between normal varieties.  

\begin{defn}[Exceptional divisors]\label{p-def2.7}
Let $f\colon X\to Y$ be a proper surjective morphism 
between normal varieties. 
Let $E$ be a Weil divisor on $X$. 
We say that $E$ is {\em{$f$-exceptional}} 
if $\mathrm{codim}_Yf(\Supp E)\geq 2$. 
We note that $f$ is not always 
assumed to be birational. 
\end{defn}

\section{Preliminary lemmas}\label{p-sec3}

In this section, we collect some useful and important lemmas for the 
reader's convenience. 
They are more or less well known to the experts. 

\medskip 

Let us start with the following 
easy lemmas on $\mathbb R$-divisors. 
We will use them repeatedly in this paper. 

\begin{lem}\label{p-lem3.1}
Let $A$ be a Cartier divisor on a normal 
variety $V$. 
Let $B$ be an $\mathbb R$-Cartier $\mathbb R$-divisor 
on $V$ such that 
$B=\sum _{i\in I} b_iB_i$ where $b_i\in \mathbb R$ and $B_i$ is a prime 
divisor on $V$ for every $i$ with 
$B_i\ne B_j$ for $i\ne j$. 
Assume that $A\sim _{\mathbb R} B$. 
Then we can take a $\mathbb Q$-Cartier $\mathbb Q$-divisor 
$C=\sum _{i\in I} c_i B_i$ on $V$ such that 
\begin{itemize}
\item[(i)] $A\sim _{\mathbb Q} C$, 
\item[(ii)] $c_i=b_i$ holds if $b_i\in \mathbb Q$, and 
\item[(iii)] $|c_i-b_i|\ll 1$ holds for $b_i\in \mathbb R\setminus 
\mathbb Q$. 
\end{itemize}
In particular, $\Supp C=\Supp B$, $\lfloor C\rfloor 
=\lfloor B\rfloor$, 
$\lceil C\rceil =\lceil B\rceil$, and $\Supp \{C\}=\Supp 
\{B\}$. 
\end{lem}

\begin{proof}
It is an easy exercise. 
For the details, see, for example, the proof of 
\cite[Lemma 4.15]{fujino-fundamental}. 
\end{proof}

\begin{lem}\label{p-lem3.2}
Let $D=\sum _{i\in I} a_i D_i$ be an $\mathbb R$-divisor 
on a smooth projective variety $V$, 
where $a_i\in \mathbb R$ and $D_i$ is a prime divisor 
on $V$ for every $i$ with $D_i\ne D_j$ for $i\ne j$. 
Assume that $D$ is semi-ample. 
Then we can construct a $\mathbb Q$-divisor 
$D^\dag=\sum _{i\in I} b_iD_i$ such that 
\begin{itemize}
\item[(i)] $D^\dag$ is semi-ample, 
\item[(ii)] $b_i=a_i$ holds if $a_i \in \mathbb Q$, and 
\item[(iii)] $|b_i-a_i|\ll1$ holds for $a_i\in \mathbb R\setminus \mathbb Q$. 
\end{itemize}
\end{lem}
\begin{proof}
Since $D$ is semi-ample, 
we can write $D=\sum _{j\in J} m_jM_j$ where $m_j$ is a 
positive real number and 
$M_j$ is a semi-ample Cartier divisor on $V$ for every $j$. As 
usual, by perturbing $m_j$s suitably, 
we get a desired semi-ample $\mathbb Q$-divisor 
$D^\dag$ on $V$. 
For the details, see, for example, 
the proof of \cite[Lemma 4.15]{fujino-fundamental}. 
\end{proof}

Next, we treat a very useful covering trick, which is essentially 
due to Yujiro Kawamata. 
We will use it in the proof of Theorem \ref{p-thm9.3}. 

\begin{lem}\label{p-lem3.3}
Let $f\colon V\to W$ be a projective surjective morphism 
between smooth quasi-projective varieties and let 
$H$ be a Cartier divisor on $W$. 
Let $d$ be an arbitrary positive integer. 
Then we can take a finite flat morphism 
$\tau\colon W'\to W$ from a smooth 
quasi-projective variety $W'$ and a Cartier divisor $H'$ on $W'$ such that 
$\tau^*H\sim dH'$ and that 
$V'=V\times _W W'$ is a smooth quasi-projective variety with 
$\omega_{V'/W'}=\rho^*\omega_{V/W}$, where $\rho\colon V'\to V$. 
By construction, we may assume that 
$\tau\colon W'\to W$ is Galois. 
\end{lem}

\begin{proof}
We take general very ample Cartier divisors $D_1$ and $D_2$ with the following 
properties. 
\begin{itemize}
\item[(i)] $H\sim D_1-D_2$, 
\item[(ii)] $D_1$, $D_2$, $f^*D_1$, and $f^*D_2$ are smooth, 
\item[(iii)] $D_1$ and $D_2$ have no common components, and 
\item[(iv)] $\Supp (D_1+D_2)$ and $\Supp (f^*D_1+f^*D_2)$ 
are simple normal crossing divisors. 
\end{itemize}
We take a finite flat cover due to Kawamata 
with respect to $W$ and $D_1+D_2$ 
(see \cite[3.19.~Lemma]{esnault-viehweg} and 
\cite[Lemma 2.5]{viehweg3}). 
Then we obtain $\tau\colon W'\to W$ and $H'$ such that $\tau^*H\sim dH'$. 
By the construction of the above Kawamata cover $\tau\colon W'\to W$, 
we may assume that the ramification locus $\Sigma$ of $\tau$ in $W$ is 
a general simple normal crossing divisor. This means that 
$f^*P$ is a smooth divisor for any irreducible component 
$P$ of 
$\Sigma$ and that $f^*\Sigma$ is a simple normal crossing divisor on $V$. 
In this situation, we can easily check that 
$V'=V\times _WW'$ is a smooth quasi-projective variety. 
$$
\xymatrix{
 V' \ar[r]^{\rho} \ar[d]_{f'} & V\ar[d]^{f} \\
  W' \ar[r]_{\tau} & W
} 
$$
By construction, we can also easily check that 
$\omega_{V'/W'}=\rho^*\omega_{V/W}$ by the Hurwitz formula. 

Let us see the construction of $f'\colon V'\to W'$ more precisely 
for the reader's convenience. 
Let $\mathcal A$ be an ample invertible sheaf on $W$ 
such that $\mathcal A^{\otimes d}\otimes \mathcal O_W(-D_i)$ is generated by 
global sections for $i=1, 2$. 
We put $n=\dim W$. 
We take smooth divisors 
$$H^{(1)}_1, \ldots, H^{(1)}_n, H^{(2)}_1, 
\ldots, H^{(2)}_n$$ on $W$ 
in general position such that 
$\mathcal A^{\otimes d}= \mathcal O_W(D_i+H^{(i)}_j)$ for 
$1\leq j\leq n$ and $i=1, 2$. 
Let $Z^{(i)}_j$ be the cyclic cover associated to 
$\mathcal A^{\otimes d}= \mathcal O_W(D_i+H^{(i)}_j)$ for 
$1\leq j\leq n$ and $i=1, 2$. 
Then $W'$ is the normalization of 
$$
\left(Z^{(1)}_1\times _W \cdots \times _WZ^{(1)}_n\right)
\times _W \left(Z^{(2)}_1\times _W \cdots \times _WZ^{(2)}_n\right). 
$$ 
We note that $W'$ is smooth since 
$$\bigcap_{j=1}^n (D_i\cap H^{(i)}_j)=\emptyset$$ for 
$i=1, 2$. 
For the details, see, for example, 
\cite[3.19.~Lemma]{esnault-viehweg} and 
\cite[Lemma 2.5]{viehweg3}. 
Let $S^{(i)}_j$ be the cyclic cover of $V$ associated 
to $\left(f^*\mathcal A\right)^{\otimes d} 
=\mathcal O_V(f^*D_i+f^*H^{(i)}_j)$. 
Then we define $V'$ as 
the normalization of 
$$
\left(S^{(1)}_1\times _V \cdots \times _VS^{(1)}_n\right)
\times _V \left(S^{(2)}_1\times _V \cdots \times _VS^{(2)}_n\right). 
$$ 
As before, $V'$ is smooth 
since 
$$
\bigcap _{j=1}^n(f^*D_i\cap f^*H^{(i)}_j)=\emptyset
$$ 
for $i=1, 2$. 
Note that $\rho\colon V'\to V$ is a finite flat morphism between smooth 
quasi-projective varieties. 
Since $V\times _WW'\to V$ is finite and flat 
and $V$ is smooth, 
$V\times _W W'$ is Cohen--Macaulay 
(see, for example, \cite[Corollary 5.5]{kollar-mori}). 
We put 
$$
\Sigma =D_1+D_2+\sum _{i, j} H^{(i)}_j. 
$$ 
Then $\Sigma$ is a simple normal crossing divisor on $W$. 
Let $U$ be the largest Zariski open set of $W$ such that 
$f$ is smooth over $U$. 
We put $$\Pi=\Sigma_{\mathrm{sing}}\cup \left((W\setminus U)
\cap \Sigma\right), 
$$ 
where $\Sigma_{\mathrm{sing}}$ is the singular locus of $\Sigma$. 
Since 
$$
D_1, D_2, H^{(1)}_1, \ldots, H^{(1)}_n, H^{(2)}_1, 
\ldots, H^{(2)}_n
$$ 
are general, $\mathrm{codim}_W\Pi\geq 2$ and 
$\mathrm{codim}_Vf^{-1}\Pi\geq 2$ hold. 
We note that $\tau\colon W'\to W$ is \'etale outside $\Sigma$. 
Let $P$ be any closed point of $\Sigma\setminus \Pi$. 
Then $f\colon V\to W$ is smooth over a neighborhood of $P$. 
Hence we can check  
that $V\times _W W'$ is smooth 
in codimension one. 
Therefore, $V\times _W W'$ is normal. 
Since $\rho$ factors through 
$V\times _W W'$ by construction, 
we see that $V'=V\times _W W'$ by Zariski's main theorem. 
By the above description, 
if $K_{W'}=\tau^*K_W+R$, 
then $K_{V'}=\rho^*K_V+f'^*R$. 
Therefore, $\omega_{V'/W'}=\rho^*\omega_{V/W}$ holds. 
By the above construction of $\tau\colon W'\to W$, 
we see that $\tau\colon W'\to W$ is Galois. 
\end{proof}

We give a very important remark on Lemma \ref{p-lem3.3}. 

\begin{rem}\label{p-rem3.4}
In the proof of Lemma \ref{p-lem3.3}, 
let $S$ be any simple normal crossing divisor on $V$. 
Then we can choose the ramification locus $\Sigma$ of $\tau$ such that 
$f^*P\not\subset S$ for any irreducible component $P$ of $\Sigma$ and that 
$f^*\Sigma\cup S$ is a simple normal 
crossing divisor on $V$. 
If we choose $\Sigma$ as above, then we obtain that 
$\rho^*S$ is a simple normal crossing divisor on $V'$. 
\end{rem}

Lemma \ref{p-lem3.5} is an elementary property of 
semi-ample $\mathbb R$-divisors. 
We give a proof for the sake of completeness. 

\begin{lem}\label{p-lem3.5} 
Let $f\colon  V\to W$ be a surjective morphism between normal 
projective varieties. 
Let $D$ be a nef and $f$-semi-ample $\mathbb R$-divisor 
on $V$ and let $H$ be an ample $\mathbb R$-divisor 
on $W$. 
Then $aD+bf^*H$ is semi-ample for any positive real numbers $a$ and 
$b$. 
\end{lem}
\begin{proof} 
Since $D$ is $f$-semi-ample, 
there exists a surjective morphism $g\colon V\to Z$ with 
the following commutative diagram 
$$
\xymatrix{
V \ar[dr]_-f\ar[rr]^-g& &Z \ar[dl]^-h\\
&W&
}
$$ 
such that 
\begin{itemize}
\item[(i)] $Z$ is a normal projective variety, and 
\item[(ii)] $D\sim _{\mathbb R} g^*A$, where $A$ is a nef 
and $h$-ample $\mathbb R$-divisor on $Z$. 
\end{itemize}
We can take a large positive real number $c$ such that 
$A+ch^*H$ is ample because $H$ is ample and 
$A$ is $h$-ample. 
If $b/a\geq c$ holds, 
then 
$$
aA+bh^*H=a(A+ch^*H)+(b-ac)h^*H
$$
is ample. 
If $b/a<c$ holds, 
then 
$$
aA+bh^*H=(a-b/c)A+(b/c)(A+ch^*H)
$$ 
is ample. 
Hence $aA+bh^*H$ is always ample. 
This implies that 
$$
aD+bf^*H\sim _{\mathbb R} g^*(aA+bh^*H)
$$ 
is semi-ample for any positive real numbers $a$ and $b$. 
\end{proof}

Finally, let us explain Viehweg's fiber product trick. 
We include the proof for the benefit of 
the reader. 
We will use it in the proof of Theorems \ref{q-thm1.8} and \ref{q-thm1.9} 
in Section \ref{p-sec10}. 

\begin{lem}[{see \cite[(4.9) Lemma]{mori}}]\label{p-lem3.6}
Let $V$ be a reduced Gorenstein scheme. Note that 
$V$ may be reducible. 
We consider 
$$
V'\overset{\delta}\longrightarrow V^\nu\overset{\nu}\longrightarrow V
$$ 
where $\nu\colon V^\nu\to V$ is the normalization and 
$\delta\colon V'\to V^\nu$ is a resolution of singularities. 
Then, for every positive integer $n$, we have 
\begin{equation}\label{p-3.1} 
\nu_*\mathcal O_{V^\nu}(nK_{V^\nu})\subset \omega^{\otimes n}_V
\end{equation}
and 
\begin{equation}\label{p-3.2}
\delta_*\mathcal O_{V'}(nK_{V'}+E)\subset \mathcal O_{V^\nu}(nK_{V^\nu})
\end{equation}
where $E$ is any $\delta$-exceptional 
divisor on $V'$. 
In particular, we have 
\begin{equation}\label{p-3.3}
(\nu\circ \delta)_*\mathcal O_{V'}(nK_{V'}+E)\subset 
\omega^{\otimes n}_V
\end{equation} 
for every positive integer $n$. 
If $U$ is a Zariski open set of $V$ such that 
$\nu\circ \delta$ is an isomorphism over 
$U$, then the inclusion 
\eqref{p-3.3} 
is an isomorphism over $U$. 
\end{lem}

\begin{proof}
In Steps \ref{p-step3.6.1} and \ref{p-step3.6.2}, 
we will prove 
\eqref{p-3.2} and \eqref{p-3.1}, respectively. 
\setcounter{step}{0}
\begin{step}\label{p-step3.6.1}
By taking the double dual of $\delta_*\mathcal O_{V'}(nK_{V'}+E)$, 
we obtain $\mathcal O_{V^\nu}(nK_{V^\nu})$. Therefore, 
we have 
$$
\delta_*\mathcal O_{V'}(nK_{V'}+E)\subset \mathcal O_{V^\nu}(nK_{V^\nu}) 
$$ 
for every integer $n$. 
\end{step} 
\begin{step}\label{p-step3.6.2}
Since $\nu$ is birational, 
the trace map 
$\nu_*\mathcal O_{V^\nu}(K_{V^\nu})\to \omega_V$ is a generically 
isomorphic injection 
\begin{equation}\label{p-3.4}
\nu_*\mathcal O_{V^\nu}(K_{V^\nu})\hookrightarrow \omega_V. 
\end{equation} 
More precisely, the above trace map is an isomorphism 
over the isomorphism locus of $\nu$ and then it is injective since $\nu_*
\mathcal O_{V^\nu}(K_{V^\nu})$ is torsion-free. 
Since $\nu$ is finite, 
\begin{equation}\label{p-3.5}
\nu^*\nu_*\mathcal O_{V^\nu}(K_{V^\nu})\to \mathcal O_{V^\nu}(K_{V^\nu})
\end{equation}
is surjective. 
Since $\mathcal O_{V^\nu}(K_{V^\nu})$ 
is torsion-free, the kernel of \eqref{p-3.5} 
is the torsion part of $\nu^*\nu_*\mathcal O_{V^\nu}(K_{V^\nu})$. 
Therefore, by \eqref{p-3.4}, we get an inclusion 
\begin{equation}\label{p-3.6}
\mathcal O_{V^\nu}(K_{V^\nu})\hookrightarrow \nu^*\omega_V 
\end{equation} 
because $\nu^*\omega_V$ is torsion-free. 
Let $n$ be a positive integer with $n\geq 2$. Then 
we have 
$$
\mathcal O_{V^\nu}(nK_{V^\nu})=\mathcal O_{V^\nu}
(K_{V^\nu}+(n-1)K_{V^\nu})\hookrightarrow 
\mathcal O_{V^\nu}(K_{V^\nu})\otimes \nu^*\omega^{\otimes n-1}_V
$$ 
by \eqref{p-3.6}. Therefore, by taking $\nu_*$, we get 
$$
\nu_*\mathcal O_{V^\nu}(nK_{V^\nu})
\hookrightarrow \nu_*\mathcal O_{V^\nu}(K_{V^\nu})\otimes 
\omega^{\otimes n-1}_V\hookrightarrow \omega^{\otimes n}_V 
$$ 
by \eqref{p-3.4}. 
This is what we wanted. 
\end{step}
By the above construction of \eqref{p-3.1} and \eqref{p-3.2}, 
it is obvious that the inclusion 
$$
(\nu\circ \delta)_*\mathcal O_{V'}(nK_{V'}+E)\subset 
\omega^{\otimes n}_V
$$ 
is an isomorphism over $U$. 
\end{proof}

\begin{lem}\label{p-lem3.7}
Let $f\colon X_0\to Y_0$ be a projective surjective morphism 
between smooth quasi-projective 
varieties and let $\Delta_0$ be an effective $\mathbb R$-divisor 
on $X_0$ such that $\Supp \Delta_0$ is a simple 
normal crossing divisor on $X_0$ 
and $(X_0, \Delta_0)$ is log canonical over a nonempty Zariski 
open set of $Y_0$. 
Let $L_0$ be a Cartier divisor on $X_0$ such that 
$L_0\sim _{\mathbb R} k(K_{X_0/Y_0}+\Delta_0)$ for 
some positive integer $k$. 
Assume that $f$ is flat. 
We consider the $s$-fold fiber product 
$$
X^s_0:=\underbrace{X_0\times _{Y_0} X_0
\times _{Y_0} \cdots \times _{Y_0} X_0}_{s}
$$ 
of $X_0$ over $Y_0$ and let $f^s\colon X^s_0\to Y_0$ be the induced 
morphism. We take a resolution of singularities $\rho\colon X^{(s)}_0\to X^s_0$ 
which 
is an isomorphism over a nonempty Zariski open set of $Y_0$. 
Then we can write 
$$
\mathcal O_{X^{(s)}_0}(K_{X^{(s)}_0})=\rho^*\omega_{X^s_0}\otimes 
\mathcal O_{X^{(s)}_0}(R)
$$ 
where $R$ is an $(f^s\circ \rho)$-vertical Cartier divisor 
by construction. 
Let $p_i\colon X^s_0\to X_0$ be the $i$-th projection. We put 
$\pi_i=p_i\circ \rho\colon X^{(s)}_0\to X_0$. 
We consider 
$$
L^{(s)}_0:=\sum _{i=1}^s \pi_i^*L_0+kR. 
$$ 
We further assume that $f_*\mathcal O_{X_0}(L_0)$ is locally free. 
Then there exists a generically isomorphic injection 
$$
f^{(s)}_*\mathcal O_{X^{(s)}_0}(L^{(s)}_0)\hookrightarrow 
\bigotimes _{i=1}^sf_*\mathcal O_{X_0}(L_0) 
$$ 
with $f^{(s)}=f^s\circ \rho$. 
We have 
$$
L^{(s)}_0
\sim _{\mathbb R}k\left(K_{X^{(s)}_0/Y_0}+\sum _{i=1}^s \pi_i^*\Delta_0\right).  
$$ 
Note that 
$$
\left(X^{(s)}_0, \sum _{i=1}^s \pi_i^*\Delta_0\right)
$$ 
is log canonical over a nonempty Zariski open set of $Y_0$. 
We also note that $X^{(s)}_0$ may be reducible, that is, 
$X^{(s)}_0$ may be a disjoint union of smooth varieties. 
\end{lem}

\begin{proof}
By the flat base change theorem, 
we have $$\omega_{X^s_0/Y_0}=\bigotimes _{i=1}^s p^*_i\omega_{X_0/Y_0}. 
$$ 
In particular, 
$X^s_0$ is Gorenstein. 
We note that 
\begin{equation}\label{p-3.7}
\begin{split}
L^{(s)}_0&=\sum _{i=1}^s \rho^*p^*_i (kK_{X_0/Y_0}+(L_0-kK_{X_0/Y_0}))+kR
\\& \sim kK_{X^{(s)}_0/Y_0}+\sum _{i=1}^s \pi^*_i (L_0-kK_{X_0/Y_0}). 
\end{split} 
\end{equation} 
\begin{claim} We have the following isomorphism of 
locally free sheaves:  
$$f^s_*\mathcal O_{X^s_0}\left(\sum _{i=1}^s p^*_iL_0\right)\simeq 
\bigotimes ^s f_*\mathcal O_{X_0}(L_0). $$
\end{claim}
\begin{proof}[Proof of Claim]
We use induction on $s$. 
If $s=1$, then the statement is obvious. 
So we assume that $s\geq 2$. 
We consider the following commutative diagram 
$$
\xymatrix{
X^s_0 \ar[dr]^-{f^s}\ar[d]_-{p_s}\ar[r]^-q& X^{s-1}_0 \ar[d]^-{f^{s-1}}\\
X_0 \ar[r]_-f& Y_0
}
$$
where $q=(p_1, \cdots, p_{s-1})$. 
Then we have 
\begin{equation}\label{p-3.8}
\mathcal O_{X^s_0}\left(\sum _{i=1}^s p^*_i L_0\right)\simeq 
\mathcal O_{X^s_0}(p^*_sL_0)\otimes 
q^*\mathcal O_{X^{s-1}_0}\left(\sum _{i=1}^{s-1}p^*_iL_0\right). 
\end{equation}
Therefore, we obtain 
\begin{equation*}
\begin{split}
f^s_*\mathcal O_{X^s_0}\left(\sum _{i=1}^s p^*_iL_0\right)
&\simeq  f_*{p_s}_*\left(\mathcal O_{X^s_0}(p^*_sL_0)\otimes 
q^*\mathcal O_{X^{s-1}_0}\left(\sum _{i=1}^{s-1}p^*_iL_0\right)\right)
\\
&\simeq 
f_*\left(\mathcal O_{X_0}(L_0)\otimes {p_s}_*q^*\mathcal O_{X^{s-1}_0}\left(
\sum _{i=1}^{s-1}p^*_iL_0\right)\right)
\\ 
&\simeq 
f_*\left(\mathcal O_{X_0}(L_0)\otimes f^*f^{s-1}_*\mathcal O_{X^{s-1}_0}\left(
\sum _{i=1}^{s-1}p^*_iL_0\right)\right)
\\ 
&\simeq 
f_*\left(\mathcal O_{X_0}(L_0)
\otimes f^*\left(\bigotimes^{s-1}f_*\mathcal O_{X_0}(L_0)\right)\right) 
\\ 
&\simeq f_*\mathcal O_{X_0}(L_0)\otimes \bigotimes ^{s-1}f_*\mathcal O_{X_0}(L_0)
\\ &\simeq 
\bigotimes ^s f_*\mathcal O_{X_0}(L_0). 
\end{split}
\end{equation*}
Note that the first isomorphism follows from \eqref{p-3.8}, 
the second one is due to the projection formula, 
the third one is obtained by the flat base change theorem, 
the fourth one is due to induction on $s$, and the fifth one follows 
from the projection formula. 
Hence we obtain the desired isomorphism. 
\end{proof}
Let us go back to the proof of Lemma \ref{p-lem3.7}. 
We have an inclusion 
\begin{equation*}
\begin{split}
\rho_*\mathcal O_{X^{(s)}_0}(L^{(s)}_0)&\subset 
\omega^{\otimes k}_{X^s_0/Y_0}\otimes \mathcal O_{X^s_0}\left(\sum 
_{i=1}^s p^*_i (L_0-kK_{X_0/Y_0})\right)\\ 
&\simeq \mathcal O_{X^s_0}\left(\sum_{i=1}^s p^*_i L_0\right)
\end{split}
\end{equation*}
by \eqref{p-3.7} and 
Lemma \ref{p-lem3.6}, which is an isomorphism 
over a nonempty Zariski open set of $Y_0$. 
By taking $f^s_*$, 
the Claim yields a generically isomorphic injection 
\begin{equation*}
\begin{split}
f^{(s)}_*\mathcal O_{X^{(s)}_0}(L^{(s)}_0)&
\subset f^s_*\mathcal O_{X^s_0}
\left(\sum_{i=1}^s p^*_i L_0\right)\\
&\simeq \bigotimes _{i=1}^s f_*\mathcal O_{X_0}(L_0), 
\end{split} 
\end{equation*} 
where $f^{(s)}=f^s\circ \rho\colon X^{(s)}_0\to Y_0$. 
By assumption, $L_0-kK_{X_0/Y_0}\sim_{\mathbb R} k\Delta_0$. 
Therefore, 
$$
L^{(s)}_0\sim _{\mathbb R} k\left(K_{X^{(s)}_0/Y_0}
+\sum _{i=1}^s \pi^*_i \Delta_0\right)
$$ 
by \eqref{p-3.7}. 
We can take a nonempty Zariski open set $U$ of $Y_0$ such 
that $f$ is smooth 
over $U$, $\Supp \Delta$ is relatively simple normal crossing 
over $U$, and $\rho$ is an isomorphism over $U$. 
Then we see that $\left(X^{(s)}_0, \sum _{i=1}^s \pi^*_i \Delta_0\right)$ 
is log canonical over $U$. 
\end{proof}

\section{Weakly positive sheaves and big sheaves}\label{p-sec4}

We briefly recall some basic properties of 
Viehweg's weakly positive sheaves and big sheaves. 
For the details, see \cite[Chapter 3]{fujino-iitaka}, 
\cite{viehweg1}, \cite{viehweg2}, and \cite{viehweg3}.  

\begin{defn}[Weak positivity and 
bigness]\label{p-def4.1}\index{weak positivity}
\index{big sheaf} 
Let $\mathcal F$ be a torsion-free coherent sheaf on a smooth 
quasi-projective variety $W$. 
We say that $\mathcal F$ is {\em{weakly positive}} 
if, for every positive integer $\alpha$ and every ample invertible 
sheaf $\mathcal H$, 
there exists a positive integer $\beta$ such that 
$\widehat{S}^{\alpha\beta}(\mathcal F)
\otimes \mathcal H^{\otimes \beta}$ is generically 
generated by global sections. 
We say that a nonzero torsion-free coherent sheaf 
$\mathcal F$ is {\em{big}} (in the sense of Viehweg) 
if, for every ample invertible sheaf $\mathcal H$, 
there exists a positive integer $a$ such that 
$\widehat{S}^a(\mathcal F)\otimes \mathcal H^{\otimes -1}$ is weakly positive. 
\end{defn}

\begin{rem}\label{p-rem4.2} 
In this paper and 
\cite{fujino-iitaka}, we adopt Viehweg's definition 
of weakly positive sheaves 
in \cite[Definition 3.1]{viehweg2}. Note that 
it is slightly weaker than \cite[Definition 2.1]{viehweg1}. 
The definition of weakly positive sheaves depends on the literature. 
Definition \ref{p-def4.1} seems to be sufficient and 
suitable for applications to the Iitaka 
conjecture (see \cite{fujino-iitaka}). 
\end{rem}

\begin{ex}\label{p-ex4.3} 
Let $D$ be a Cartier divisor on a smooth projective variety $W$. 
Then $\mathcal O_W(D)$ is weakly positive if and only if 
$D$ is pseudo-effective in the usual sense. 
We note that $\mathcal O_W(D)$ is big in the 
sense of Definition \ref{p-def4.1} 
if and only if $D$ is big in the usual sense, that is, 
$\kappa (W, D)=\dim W$. 
\end{ex}

For the reader's convenience, let us recall the 
following basic properties of big sheaves without proof. 

\begin{lem}[{\cite[Lemma 3.6]{viehweg2} and \cite[Lemma 3.1.15]{fujino-iitaka}}]\label{p-lem4.4}
Let $\mathcal F$ be a nonzero torsion-free coherent sheaf on a smooth 
quasi-projective variety $W$. 
Then the following three conditions are equivalent. 
\begin{itemize}
\item[(i)] There exist an ample invertible sheaf $\mathcal H$ on $W$, 
some positive integer $\nu$, and an inclusion $\bigoplus \mathcal H\hookrightarrow 
\widehat{S}^\nu(\mathcal F)$, 
which is an isomorphism over a nonempty Zariski open set of $W$. 
\item[(ii)] For every invertible sheaf $\mathcal M$ on $W$, 
there exists some positive integer $\gamma$ such that 
$\widehat {S}^\gamma (\mathcal F)\otimes \mathcal M^{\otimes -1}$ is weakly 
positive. In particular, $\mathcal F$ is a big sheaf.  
\item[(iii)] There exist some positive integer $\gamma$ and 
an ample invertible sheaf $\mathcal M$ on $W$ such that 
$\widehat {S}^\gamma (\mathcal F)\otimes \mathcal M^{\otimes -1}$ is weakly 
positive. 
\end{itemize}
\end{lem}

We will use the following three easy lemmas on big sheaves 
in this paper. 
So we explicitly state them here for the reader's convenience. 

\begin{lem}\label{p-lem4.5}
Let $\mathcal F$ be a weakly positive sheaf and let $\mathcal H$ be an 
ample invertible sheaf on a smooth quasi-projective 
variety $W$. 
Then $\mathcal F\otimes \mathcal H$ is big. 
\end{lem}

We give a proof for the sake of completeness. 

\begin{proof}[Proof of Lemma \ref{p-lem4.5}]
Since $\mathcal F$ is weakly positive, 
$\widehat{S}^{2b}(\mathcal F)\otimes \mathcal H^{\otimes b}$ 
is generically generated by 
global sections for some positive integer $b$. 
By replacing $b$ with a multiple, 
we may assume that 
$\mathcal H^{\otimes b-1}$ 
is generated by global sections. 
Then $\widehat{S}^{2b}(\mathcal F\otimes \mathcal H)
\otimes \mathcal H^{\otimes -1}=
\widehat{S}^{2b}(\mathcal F)\otimes \mathcal H^{\otimes 2b-1}$ 
is generically generated by global sections. 
In particular, $\widehat {S}^{2b}(\mathcal F\otimes \mathcal H)
\otimes \mathcal H^{\otimes -1}$ is weakly positive. 
This implies that 
$\mathcal F\otimes \mathcal H$ is big by 
Lemma \ref{p-lem4.4}. 
\end{proof}

\begin{lem}\label{p-lem4.6} 
Let $\mathcal G$ be any nonzero torsion-free 
coherent sheaf on a smooth 
quasi-projective variety $W$ and let $\mathcal H$ be 
an ample invertible sheaf on $W$. 
Then there exists a positive integer $l$ 
such that $\mathcal G\otimes \mathcal H^{\otimes l}$ 
is big.  
\end{lem}
\begin{proof} 
Since $\mathcal H$ is ample, 
$\mathcal G\otimes \mathcal H^{\otimes m}$ is 
generated by global sections for some 
positive integer $m$. 
Therefore, there exists a surjection 
$$
\bigoplus _{\mathrm{finite}}\mathcal O_W\to \mathcal 
G \otimes \mathcal H^{\otimes m}. 
$$ 
This implies that $\mathcal G\otimes \mathcal H^{\otimes m}$ 
is weakly positive (see \cite[Lemma 3.1.12 (ii)]{fujino-iitaka}). 
By Lemma \ref{p-lem4.5}, we obtain that 
$\mathcal G\otimes \mathcal H^{\otimes l}$ 
is big for every integer $l\geq m+1$. 
\end{proof}

\begin{lem}\label{p-lem4.7}
Let $\mathcal F$ be a torsion-free coherent sheaf on a smooth 
quasi-projective variety $W$ and let $\tau\colon W'\to W$ be 
a finite surjective morphism from a smooth quasi-projective 
variety $W'$. 
Assume that $\tau^*\mathcal F$ is big. 
Then $\mathcal F$ is a big sheaf on $W$. 
\end{lem}

We include the proof for the benefit of the reader. 

\begin{proof}[Proof of Lemma \ref{p-lem4.7}] 
We take an ample invertible sheaf $\mathcal H$ on $W$. 
By replacing $W$ with $W\setminus \Sigma$ for 
some suitable closed subset $\Sigma$ of 
codimension $\geq 2$ (see, for example, 
\cite[Lemma 3.1.12 (i)]{fujino-iitaka}), 
we may assume that $\mathcal F$ is locally free. 
Since $\tau^*\mathcal F$ is big by assumption, 
there exists a positive integer $a$ such that 
$S^a(\tau^*\mathcal F)\otimes \tau^*\mathcal H^{\otimes -1}=
\tau^*\left(S^a(\mathcal F)\otimes 
\mathcal H^{\otimes -1}\right)$ is weakly positive 
(see Lemma \ref{p-lem4.4}). 
Therefore, $S^a(\mathcal F)\otimes 
\mathcal H^{\otimes -1}$ is weakly positive since $\tau$ is 
finite 
(see, for example, \cite[Lemma 3.1.12 (v)]{fujino-iitaka}). 
This means that $\mathcal F$ is big by 
Lemma \ref{p-lem4.4}.  
\end{proof}

\section{Mixed-$\omega$-sheaves and mixed-$\omega$-big sheaves}
\label{p-sec5} 

In this section, we introduce mixed-$\omega$-sheaves, 
mixed-$\omega$-big-sheaves, mixed-$\widehat\omega$-sheaves, 
and mixed-$\widehat\omega$-big-sheaves. 
We also treat some important examples in Lemmas \ref{p-lem5.5}, 
\ref{p-lem5.9}, and \ref{p-lem5.10}. 

\medskip 

Let us start with the definition of mixed-$\omega$-sheaves and 
pure-$\omega$-sheaves. 

\begin{defn}[Mixed-$\omega$-sheaf and 
pure-$\omega$-sheaf]\label{p-def5.1} 
A torsion-free coherent sheaf $\mathcal F$ on a normal quasi-projective 
variety $W$ is called a {\em{mixed-$\omega$-sheaf}} 
if there exist a projective surjective morphism 
$f\colon V\to W$ from a smooth quasi-projective variety $V$ and a simple 
normal crossing divisor $D$ on $V$ such that 
$\mathcal F$ is a direct summand of $f_*\mathcal O_V(K_V+D)$. 
When $D=0$, $\mathcal F$ is called a {\em{pure-$\omega$-sheaf}} on $W$. 
\end{defn}

We give a very important remark on Definition \ref{p-def5.1}. 

\begin{rem}[Pure-$\omega$-sheaves versus 
Nakayama's $\omega$-sheaves]\label{p-rem5.2}
The notion of pure-$\omega$-sheaves is essentially the same 
as that of Nakayama's $\omega$-sheaves in \cite{nakayama} 
when we treat torsion-free coherent sheaves on normal projective 
varieties 
(see the Remark after \cite[Chapter V, 3.8.~Definition]{nakayama}). 
However, the definition of pure-$\omega$-sheaves 
in Definition \ref{p-def5.1} does not coincide with 
\cite[Chapter V, 3.8.~Definition]{nakayama}. 
Our definition seems to be more reasonable than Nakayama's 
from the mixed Hodge theoretic viewpoint. 
\end{rem}

For some geometric applications, the notion of 
mixed-$\omega$-big-sheaves and pure-$\omega$-big-sheaves 
is very useful. 

\begin{defn}[Mixed-$\omega$-big-sheaf and 
pure-$\omega$-big-sheaf]\label{p-def5.3}
Let $\mathcal F$ be a torsion-free coherent sheaf 
on a normal quasi-projective 
variety $W$. 
If there exist projective 
surjective morphisms $f\colon V\to W$, 
$p\colon V\to Z$, and an ample divisor $A$ on $Z$ 
satisfying the following conditions:
\begin{itemize}
\item[(i)] $V$ is a smooth quasi-projective variety, 
\item[(ii)] $Z$ is a normal quasi-projective variety, 
\item[(iii)] $D$ is a simple normal crossing divisor on $V$, 
\item[(iv)] there exists a projective surjective morphism 
$q\colon Z\to W$ such that $f=q\circ p$, and 
\item[(v)] $\mathcal F$ is a direct summand of 
$f_*\mathcal O_V(K_V+D+P)$, where $P$ is 
a Cartier divisor on $V$ such that 
$P\sim _{\mathbb Q}p^*A$,  
\end{itemize}
then $\mathcal F$ is called a {\em{mixed-$\omega$-big-sheaf}} on $W$. 
As in Definition \ref{p-def5.1}, 
$\mathcal F$ is called a {\em{pure-$\omega$-big-sheaf}} on 
$W$ when $D=0$. 
The relationships between $V$, $W$, $Z$ and 
$f$, $p$, $q$ can be visualized as follows. 
$$
\xymatrix{
V \ar[dd]_-f\ar[dr]^-p& \\ 
& Z \ar[dl]^-q\\
W
}
$$
\end{defn}

\begin{rem}\label{p-rem5.4}
Of course, we defined mixed-$\omega$-big-sheaves and 
pure-$\omega$-big-sheaves referring to 
\cite[Chapter V, 3.16.~Definition (1)]{nakayama}. 
However, Nakayama's definition of $\omega$-bigness 
is different from ours. 
Roughly speaking, we treat only a special case 
where $X=Y$ in 
\cite[Chapter V, 3.16.~Definition (1)]{nakayama}. 
\end{rem}

Lemma \ref{p-lem5.5} gives 
a very basic example of mixed-$\omega$-sheaves. 

\begin{lem}\label{p-lem5.5}
Let $V$ be a smooth quasi-projective variety 
and let $D$ be a simple normal crossing divisor on $V$. 
Let $L$ be a semi-ample 
Cartier divisor on $V$. 
Then $\mathcal O_V(K_V+D+L)$ is a mixed-$\omega$-sheaf 
on $V$ and $\mathcal O_V(K_V+L)$ is a pure-$\omega$-sheaf on $V$.  
\end{lem}

Although this lemma is well known, we give a proof for the 
sake of completeness. 

\begin{proof}[Proof of Lemma \ref{p-lem5.5}]
Let $m$ be a positive integer such that 
$|mL|$ is free. We take a general section 
$s\in H^0(V, \mathcal O_V(mL))$, whose zero divisor is $B$. 
We may assume that $B$ is a smooth divisor, $B$ and $D$ have no common 
irreducible components, and 
$\Supp(B+D)$ is a simple normal crossing divisor 
on $V$. 
The dual of 
$$
s\colon  \mathcal O_V\to \mathcal O_V(mL)
$$ 
defines an $\mathcal O_V$-algebra structure 
on 
$$
\bigoplus_{i=0}^{m-1}\mathcal O_V(-iL). 
$$ 
We put 
$$\pi\colon Z:=\Spec _V\bigoplus_{i=0}^{m-1}\mathcal O_V(-iL)
\to V. $$ 
Then $Z$ is a smooth quasi-projective variety and $\pi^*D$ is a simple 
normal crossing divisor on $Z$ by construction. 
We can check that 
$$
\pi_*\mathcal O_Z(K_Z+\pi^*D)\simeq \bigoplus _{i=0}^{m-1}
\mathcal O_V(K_V+D+iL)
$$ 
since $\pi_*\mathcal O_V=\bigoplus_{i=0}^{m-1}
\mathcal O_V(-iL)$. 
This means that 
$\mathcal O_V(K_V+D+L)$ is a mixed-$\omega$-sheaf on $V$. 
We put $D=0$ in the above argument. 
Then we see that $\mathcal O_V(K_V+L)$ is a pure-$\omega$-sheaf on $V$. 
\end{proof}

We treat two elementary lemmas. 

\begin{lem}\label{p-lem5.6}
Let $\mathcal F$ be a 
mixed-$\omega$-big-sheaf {\em{(}}resp.~pure-$\omega$-big-sheaf{\em{)}} 
on a normal quasi-projective variety $W$. 
Then $\mathcal F$ is a 
mixed-$\omega$-sheaf {\em{(}}resp.~pure-$\omega$-sheaf{\em{)}} on $W$. 
\end{lem}

\begin{proof}
We may assume that $\mathcal F$ is a direct summand of 
$f_*\mathcal O_V(K_V+D+P)$ as in Definition 
\ref{p-def5.3}. 
By Lemma \ref{p-lem5.5}, 
$\mathcal O_V(K_V+D+P)$ is 
a mixed-$\omega$-sheaf on $V$. 
Therefore, we see that $\mathcal F$ is a mixed-$\omega$-sheaf 
on $W$. If we put $D=0$, then we see that 
$\mathcal F$ is a pure-$\omega$-sheaf on $W$. 
\end{proof}

\begin{lem}\label{p-lem5.7}
Let $\mathcal F$ be a mixed-$\omega$-sheaf 
{\em{(}}resp.~pure-$\omega$-sheaf{\em{)}} 
on a normal quasi-projective variety $W$ and 
let $\mathcal A$ be an ample invertible sheaf on $W$. 
Then $\mathcal F\otimes \mathcal A$ is 
a mixed-$\omega$-big-sheaf {\em{(}}resp.~pure-$\omega$-big-sheaf{\em{)}} on $W$. 
\end{lem}

\begin{proof}
We may assume that $\mathcal F$ 
is a direct summand of 
$f_*\mathcal O_V(K_V+D)$ as in 
Definition \ref{p-def5.1}. 
We put $Z=W$. Let $A$ be an ample 
divisor on $W$ such that 
$\mathcal O_W(A)=\mathcal A$. 
Then $\mathcal F\otimes 
\mathcal A$ is a direct summand of 
$f_*\mathcal O_V(K_V+D+f^*A)$. 
Therefore, $\mathcal F\otimes 
\mathcal A$ is a mixed-$\omega$-big-sheaf 
on $W$. When $D=0$, we see that 
$\mathcal F\otimes \mathcal A$ is a pure-$\omega$-big-sheaf on $W$. 
\end{proof}

We can not replace the assumption that $\mathcal A$ 
is ample 
with one that $\mathcal A$ is big in Lemma \ref{p-lem5.7}. 

\begin{rem}\label{p-rem5.8} 
Let $\mathcal F$ be a pure-$\omega$-sheaf on a normal 
quasi-projective variety $W$ and let $B$ be an effective 
big Cartier divisor on $W$. 
A simple example (see Example \ref{p-ex7.2} below) 
shows that $\mathcal F\otimes \mathcal O_W(B)$ 
is 
not necessarily a mixed-$\omega$-big-sheaf on $W$. 
\end{rem}

Lemmas \ref{p-lem5.9} and \ref{p-lem5.10} 
give many nontrivial important examples of 
mixed-$\omega$-sheaves and mixed-$\omega$-big-sheaves 
in the study of higher-dimensional algebraic varieties.  

\begin{lem}\label{p-lem5.9} 
Let $f\colon  V\to W$ be a projective surjective morphism from 
a smooth projective variety $V$ onto a normal projective variety $W$. 
Let $D$ be a simple normal crossing divisor on $V$ and 
let $M$ be an $\mathbb R$-divisor on $V$ such that 
$M-f^*H$ is semi-ample 
for some 
ample $\mathbb Q$-divisor $H$ on $W$. 
We assume that $D$ and $\Supp \{M\}$ have no 
common irreducible components and $\Supp(D+\{M\})$ 
is a 
simple normal crossing divisor on $V$. 
Then $f_*\mathcal O_V(K_V+D+\lceil M\rceil)$ is a mixed-$\omega$-big-sheaf 
on $W$ and $f_*\mathcal O_V(K_V+\lceil M\rceil)$ 
is a pure-$\omega$-big-sheaf on $W$. 
\end{lem}

\begin{proof}
By Lemma \ref{p-lem3.2}, 
we can construct a $\mathbb Q$-divisor $M^\dag$ 
on $V$ such that $M^\dag-f^*H$ is 
semi-ample, $\Supp \{M^\dag\}=
\Supp \{M\}$, and $\lceil M^\dag\rceil =\lceil M\rceil$. 
Therefore, we may assume that $M$ is a 
$\mathbb Q$-divisor by replacing 
$M$ with $M^\dag$. 
By Kawamata's covering construction, 
we can construct a finite Galois cover $\pi\colon V'\to V$ from a smooth 
projective variety $V'$ with the following properties: 
\begin{itemize}
\item[(i)] $\pi^*D$ is a simple normal crossing divisor on $V'$, 
\item[(ii)] $\pi^*\{M\}$ is a $\mathbb Z$-divisor 
on $V'$, 
\item[(iii)] $\Supp (\pi^*D+\pi^*\{M\})$ is a simple 
normal crossing divisor on $V'$, and 
\item[(iv)] $\left(\pi_*\mathcal O_{V'}(K_{V'}+\pi^*D+\pi^*M)\right)^G
\simeq \mathcal O_V(K_V+D+\lceil M\rceil)$, where $G$ is the Galois group 
of $\pi\colon V'\to V$. 
\end{itemize}
By assumption, $\pi^*M$ is semi-ample. 
Let us consider the contraction morphism 
$p\colon V'\to Z$ associated to 
$|m\pi^*M|$ for some sufficiently large and 
divisible positive integer $m$. 
Since $\pi^*M-(f\circ \pi)^*H$ is semi-ample, 
we have 
a morphism $q\colon Z\to W$ with 
the following commutative diagram: 
$$
\xymatrix{
V'\ar[dd]_-{f\circ\pi}\ar[dr]^-p & \\
& Z\ar[dl]^-q\\ 
W
}
$$such that 
\begin{itemize}
\item[(a)] $Z$ is a normal projective variety, and 
\item[(b)] there is an ample $\mathbb Q$-divisor 
$A$ on $Z$ with 
$\pi^*M\sim _{\mathbb Q}p^*A$. 
\end{itemize}
Therefore, 
$f_*\mathcal O_V(K_V+D+\lceil M\rceil)$ 
is a mixed-$\omega$-big-sheaf on $W$ since 
it is a direct summand of $(f\circ \pi)_*\mathcal O_{V'}(K_{V'}+
\pi^*D+\pi^*M)$. 
We put $D=0$ in the above argument. Then 
$f_*\mathcal O_V(K_V+\lceil M\rceil)$ is a 
pure-$\omega$-big-sheaf on $W$. 
\end{proof}

\begin{lem}\label{p-lem5.10} 
Let $V$ be a smooth 
quasi-projective variety and let $D$ be a simple 
normal crossing divisor on $V$. 
Let $B$ be a $\mathbb Q$-divisor 
on $V$ such that 
$rB\sim 0$ for some positive integer $r$, 
$\Supp \{B\}$ and $D$ have no common irreducible 
components, and 
$\Supp (\{B\}+D)$ 
is a simple normal crossing divisor on $V$. 
Then there exist a generically finite proper morphism 
$\pi\colon V'\to V$ from a smooth 
quasi-projective variety $V'$ and a simple 
normal crossing divisor $D'$ on $V'$ such that 
$\mathcal O_V(K_V+D+\lceil B\rceil)$ is a 
direct summand of $\pi_*\mathcal O_{V'}(K_{V'}+D')$. 
In particular, $\mathcal O_V(K_V+D+\lceil B\rceil)$ 
is a mixed-$\omega$-sheaf on $V$. 
When $D=0$, 
$\mathcal O_V(K_V+\lceil B\rceil)$ is obviously 
a pure-$\omega$-sheaf on $V$. 
\end{lem}

\begin{proof}
If $B\sim 0$, then there is nothing to prove. 
By replacing $r$ suitably, we may assume that 
$iB\not\sim 0$ for $1\leq i\leq r-1$ and that 
$r\geq 2$. 
We consider the following 
$\mathcal O_V$-algebra 
$$
\mathcal A=\bigoplus _{i=0}^{r-1} \mathcal O_V(\lfloor 
-iB\rfloor)
$$ 
defined by an isomorphism 
$\mathcal O_V(-rB)\simeq \mathcal O_V$. 
Let $Z$ be the normalization 
of $\Spec _V\mathcal A$. Then we have 
$$
\tau_*\mathcal O_Z(K_Z+\tau^*D)\simeq \bigoplus _{i=0}^{r-1}
\mathcal O_V(K_V+D+
\lceil iB\rceil)
$$ 
where $\tau\colon Z\to V$. 
By construction, we see that 
$(Z, \tau^*D)$ is dlt. 
We take a suitable resolution of singularities $\rho\colon V'\to Z$ and 
write 
$$
K_{V'}+D'=\rho^*(K_Z+\tau^*D)+E
$$
where $D'$ is a reduced simple normal crossing divisor 
on $V'$ and 
$E$ is an effective $\rho$-exceptional 
$\mathbb Q$-divisor on $V'$. 
We put $\pi:=\tau\circ \rho\colon  V'\to V$. 
Then 
\begin{equation*}
\begin{split}
\pi_*\mathcal O_{V'}(K_{V'}+D')&\simeq 
\tau_*\mathcal O_Z(K_Z+\tau^*D)\\ 
&\simeq \bigoplus_{i=0}^{r-1}\mathcal O_V(K_V+D+\lceil iB\rceil).  
\end{split} 
\end{equation*}
Therefore, we have the desired statement. 
\end{proof}

We close this section with the definition of 
mixed-$\widehat\omega$-sheaves, 
mixed-$\widehat\omega$-big-sheaves, pure-$\widehat\omega$-sheaves, 
and pure-$\widehat\omega$-big-sheaves. 

\begin{defn}[Mixed-$\widehat\omega$-sheaf, 
mixed-$\widehat\omega$-big-sheaf, pure-$\widehat\omega$-sheaf, 
and pure-$\widehat\omega$-big-sheaf]\label{p-def5.11}
A torsion-free coherent sheaf $\mathcal G$ on a normal 
quasi-projective variety $W$ is called 
a {\em{mixed-$\widehat\omega$-sheaf}} 
(resp.~{\em{mixed-$\widehat\omega$-big-sheaf}}) 
if there exist 
a mixed-$\omega$-sheaf 
(resp.~mixed-$\omega$-big-sheaf) $\mathcal F$ on $W$ 
and a generically 
isomorphic injection $\mathcal F\hookrightarrow \mathcal G^{**}$ 
into the 
double dual $\mathcal G^{**}$ of $\mathcal G$. 
If $\mathcal F$ is a pure-$\omega$-sheaf (resp.~pure-$\omega$-big-sheaf) 
in the above inclusion $\mathcal F\hookrightarrow \mathcal G^{**}$, 
then $\mathcal G$ is called a {\em{pure-$\widehat\omega$-sheaf}} 
(resp.~{\em{pure-$\widehat\omega$-big-sheaf}}). 
\end{defn}

\begin{rem}\label{p-rem5.12} 
Let $X$ be a smooth projective variety, 
let $D$ be a simple normal crossing divisor on $X$, 
let $H$ be an ample Cartier divisor on $X$, 
and let $B$ be an effective Cartier divisor on $X$. 
Then $\mathcal O_X(K_X)$ is a pure-$\omega$-sheaf, 
$\mathcal O_X(K_X+H)$ is a pure-$\omega$-big-sheaf, 
$\mathcal O_X(K_X+B)$ is a pure-$\widehat\omega$-sheaf, 
and $\mathcal O_X(K_X+H+B)$ is a pure-$\widehat\omega$-big-sheaf. 
By definition, it is obvious that 
$\mathcal O_X(K_X+D)$ is a mixed-$\omega$-sheaf, 
$\mathcal O_X(K_X+D+H)$ is a mixed-$\omega$-big-sheaf, 
$\mathcal O_X(K_X+D+B)$ is a mixed-$\widehat\omega$-sheaf, 
and $\mathcal O_X(K_X+D+H+B)$ is a 
mixed-$\widehat\omega$-big-sheaf. 
\end{rem}

Let $f\colon X\to Y$ be a surjective morphism 
between smooth 
projective varieties and let $\Delta$ be a simple normal crossing 
divisor on $X$. 
Let $k$ be a positive integer with $k\geq 2$ and let 
$H$ be an ample Cartier divisor on $Y$. 
Then we will show that 
$$
\mathcal O_Y(K_Y+H)\otimes f_*\mathcal O_X(k(K_{X/Y}+\Delta))
$$
is a mixed-$\widehat\omega$-big-sheaf on $Y$ 
when $f_*\mathcal O_X(k(K_{X/Y}+\Delta))\ne 0$. 
This is a special case of Theorem \ref{p-thm9.3}, which 
we call a fundamental theorem of the theory of 
mixed-$\omega$-sheaves. 

\section{Basic properties:~Part 1}\label{p-sec6}

In this section, we treat the weak positivity and the bigness 
of mixed-$\omega$-sheaves and mixed-$\omega$-big-sheaves, 
respectively. 

\medskip 

Let us start with the following weak positivity theorem, 
which follows from the theory of mixed Hodge structures. 

\begin{thm}\label{p-thm6.1}
Let $f\colon V\to W$ be a projective surjective morphism 
between smooth quasi-projective varieties. 
Let $D$ be a simple normal crossing divisor on $V$. 
Then $f_*\mathcal O_V(K_{V/W}+D)$ is weakly positive.  
\end{thm}

\begin{proof}
We may assume that 
$V$ and $W$ are smooth projective varieties by 
compactifying $f\colon V\to W$ suitably. 
Then this result is more or less well known. 
For the proof based on the theory of variations of mixed Hodge structure 
(see \cite{fujino-higher}, \cite{ffs}, \cite{fujino-fujisawa}, \cite{fujisawa}, 
and so on), 
see \cite[Theorem 7.8 and Corollary 7.11]{fujino-weak}. 
For the proof based on the vanishing theorem, 
see \cite[Theorem 8.4]{fujino-quasi-alb}. 
\end{proof}

As an easy consequence of Theorem \ref{p-thm6.1}, we have: 

\begin{thm}[Weak positivity]\label{p-thm6.2}
Let $\mathcal F$ be a mixed-$\omega$-sheaf 
on a smooth quasi-projective 
variety $W$. 
Then $\mathcal F\otimes \omega^{\otimes -1}_W$ is 
weakly positive. 
\end{thm}

\begin{proof}
We may assume that $\mathcal F$ is 
a direct summand of $f_*\mathcal O_V(K_V+D)$ as 
in Definition \ref{p-def5.1}. 
By Theorem \ref{p-thm6.1}, 
$f_*\mathcal O_V(K_{V/W}+D)$ is weakly positive. 
Then $\mathcal F\otimes \omega^{\otimes -1}_W$ 
is weakly positive since it is a direct summand of $f_*\mathcal O_V(K_{V/W}+D)$. 
\end{proof}

When $\mathcal F$ is a mixed-$\omega$-big-sheaf on $W$ in 
Theorem \ref{p-thm6.2},  
$\mathcal F\otimes \omega^{\otimes -1}_W$ is not only 
weakly positive but also big.  

\begin{thm}[Bigness]\label{p-thm6.3}
Let $\mathcal F$ be a mixed-$\omega$-big-sheaf 
on a smooth quasi-projective variety $W$. 
Then $\mathcal F\otimes \omega^{\otimes -1}_W$ is big. 
\end{thm}

\begin{proof}Without loss of generality, we may assume that 
$\mathcal F$ is a direct summand of 
$f_*\mathcal O_V(K_V+D+P)$ as in Definition \ref{p-def5.3}. 
It is sufficient to prove that $f_*\mathcal O_V(K_{V/W}+D+P)$ 
is big. Let 
$$
\xymatrix{
V \ar[dd]_-f\ar[dr]^-p& \\ 
& Z \ar[dl]^-q\\
W
}
$$ and $A$ be as in Definition \ref{p-def5.3}. 
Let $H$ be an ample Cartier divisor 
on $W$. We take a positive integer $m$ such that 
$mA-q^*H$ is ample. 
We can take a finite surjective morphism $\tau\colon W'\to W$ from 
a smooth quasi-projective variety $W'$ and get the following 
commutative diagram 
$$
\xymatrix{
V' \ar[d]_-{f'}\ar[r]^-\rho& V \ar[d]^-f\\
W' \ar[r]_-\tau&W
}
$$ 
such that 
$\tau^*H\sim mH'$ for some Cartier divisor $H'$, $V'=V\times _W W'$ 
is a smooth quasi-projective variety, $\rho^*D$ is a simple 
normal crossing divisor, and $\rho^*\omega^{\otimes n}_{V/W}
=\omega^{\otimes n}_{V'/W'}$ holds for every integer $n$ 
(see Lemma \ref{p-lem3.3} and Remark \ref{p-rem3.4}). 
By Lemma \ref{p-lem4.7}, It is sufficient to 
prove that 
$$
\tau^*f_*\mathcal O_V(K_{V/W}+D+P)\simeq 
f'_*\mathcal O_{V'}(K_{V'/W'}+
\rho^*D+\rho^*P)
$$ 
is a big sheaf on $W'$. 
By construction, we see that 
$\rho^*P-f'^*H'$ is a semi-ample Cartier divisor 
on $V'$ since it is $\mathbb Q$-linearly equivalent to 
$\rho^*p^*(A-(1/m)q^*H)$. 
Therefore, by Lemma \ref{p-lem5.5}, 
$\mathcal O_{V'}(K_{V'}+\rho^*D+\rho^*P-f'^*H')$ is a mixed-$\omega$-sheaf 
on $V'$. 
Thus, $\mathcal E:=f'_*\mathcal O_{V'}(K_{V'}+\rho^*D
+\rho^*P-f'^*H')$ is 
a mixed-$\omega$-sheaf on $W'$. 
We note that 
$$
f'_*\mathcal O_{V'}(K_{V'/W'}+\rho^*D
+\rho^*P)\simeq \mathcal E \otimes \omega^{\otimes -1}_{W'}\otimes 
\mathcal O_{W'}(H'). 
$$ 
By Theorem \ref{p-thm6.2}, 
$\mathcal E\otimes \omega^{\otimes -1}_{W'}$ is weakly positive. 
By Lemma \ref{p-lem4.5}, 
$\mathcal E\otimes 
\omega^{\otimes -1}_{W'}\otimes 
\mathcal O_{W'}(H')$ is 
big since $H'$ is ample. 
This is what we wanted. 
\end{proof}

We close this section with an obvious corollary. 

\begin{cor}\label{p-cor6.4}
Let $\mathcal F$ be a mixed-$\widehat\omega$-sheaf 
{\em{(}}resp.~mixed-$\widehat\omega$-big-sheaf{\em{)}} 
on a smooth quasi-projective variety $W$. 
Then $\mathcal F\otimes \omega^{\otimes -1}_W$ is weakly 
positive {\em{(}}resp.~big{\em{)}}. 
\end{cor}

\begin{proof}
We note that $\mathcal F\otimes \omega^{\otimes -1}_W$ is 
weakly positive (resp.~big) if and 
only if so is $\mathcal F^{**}\otimes \omega^{\otimes -1}_W$. 
Therefore, the desired statement follows from Theorems \ref{p-thm6.2} and 
\ref{p-thm6.3}. 
\end{proof}

\section{Basic properties:~Part 2}\label{p-sec7}

In this section, we discuss some vanishing theorems 
for mixed-$\omega$-sheaves and several related topics. 

\begin{lem}[Vanishing theorem 
for mixed-$\omega$-big-sheaf]
\label{p-lem7.1}
Let $\mathcal F$ be a mixed-$\omega$-big-sheaf 
on a normal projective variety $W$. 
Then $H^i(W, \mathcal F\otimes \mathcal N)=0$ for 
every $i>0$ and every nef invertible sheaf $\mathcal N$ 
on $W$.  
\end{lem}

\begin{proof}
We may assume that $\mathcal F$ is a direct 
summand of $f_*\mathcal O_V(K_V+D+P)$ as in 
Definition \ref{p-def5.3}. 
Let $N$ be a Cartier divisor on $W$ such that 
$\mathcal N\simeq \mathcal O_W(N)$. 
It is sufficient to prove that 
$H^i(W, f_*\mathcal O_V(K_V+D+P+f^*N))=0$ for 
every $i>0$. 
We take an ample $\mathbb Q$-divisor 
$H$ on $W$ such that 
$A-q^*H$ is an ample 
$\mathbb Q$-divisor on $Z$, where $A$ and $q\colon Z\to W$ are as in 
Definition \ref{p-def5.3}. 
Then we can take a boundary $\mathbb Q$-divisor 
$\Delta$ on $V$ such that 
$\Delta\sim _{\mathbb Q} D+P-f^*H$ and 
that $\Supp \Delta$ is a simple normal crossing 
divisor on $V$. 
Then we have 
$$
K_V+D+P+f^*N-(K_V+\Delta)\sim _{\mathbb Q} f^*(H+N). 
$$ 
We note that $H+N$ is ample. 
Therefore, by \cite[Theorem 6.3 (ii)]{fujino-fundamental} 
(see also \cite[Theorem 3.16.3 (ii) and 
Theorem 5.6.2 (ii)]{fujino-foundations}, and so on), 
we obtain that 
$H^i(W, f_*\mathcal O_V(K_V+D+P+f^*N))=0$ for 
every $i>0$. 
\end{proof}

Example \ref{p-ex7.2} shows that the vanishing theorem 
does not necessarily 
hold for mixed-$\widehat \omega$-big-sheaves. 

\begin{ex}\label{p-ex7.2} 
We put $X=\mathbb P_{\mathbb P^1}
(\mathcal O_{\mathbb P^1}\oplus 
\mathcal O_{\mathbb P^1}(1))$. 
Let $C$ be the unique $(-1)$-curve on 
$X$. 
It is not difficult to see that there exists 
an ample Cartier divisor $H$ on $X$ such that $C\cdot H=1$. 
By definition, $\mathcal O_X(K_X+C)$ is a 
mixed-$\omega$-sheaf on $X$. 
Then, by Lemma \ref{p-lem5.7}, 
$\mathcal O_X(K_X+C+H)$ is a 
mixed-$\omega$-big-sheaf on $X$. 
By definition, 
$\mathcal O_X(K_X+C+H+C)$ is a 
mixed-$\widehat \omega$-big-sheaf on $X$. 
Let us consider the following short exact sequence 
$$
0\to \mathcal O_X(K_X+C+H)\to 
\mathcal O_X(K_X+C+H+C)\to 
\mathcal O_C(K_C+H|_C+C|_C)\to 0. 
$$
Since $C\simeq \mathbb P^1$ and $H\cdot C=1$, 
we obtain 
$$
\mathcal O_C(K_C+H|_C+C|_C)\simeq \mathcal O_
{\mathbb P^1}(K_{\mathbb P^1}). 
$$ 
Since $H^i(X, \mathcal O_X(K_X+C+H))=0$ for 
$i=1$ and $2$, 
we have 
$$
H^1(X, \mathcal O_X(K_X+C+H+C))\simeq 
H^1(\mathbb P^1, \mathcal O_{\mathbb P^1}
(K_{\mathbb P^1}))=\mathbb C. 
$$
In particular, $\mathcal O_X(K_X+C+H+C)$ is not a 
mixed-$\omega$-big-sheaf on $X$ by Lemma \ref{p-lem7.1}. 
We note that $\mathcal O_X(K_X)$ is a pure-$\omega$-sheaf 
and that $H+2C$ is an effective big Cartier divisor. 
However, $\mathcal O_X(K_X+H+2C)$ is not a 
mixed-$\omega$-big-sheaf. 
\end{ex}

As an easy consequence of Lemma \ref{p-lem7.1}, we have: 

\begin{lem}\label{p-lem7.3}
Let $\mathcal F$ be a mixed-$\omega$-sheaf {\em{(}}resp.~mixed-$\omega$-big-sheaf{\em{)}} on a normal projective variety $W$ with $\dim W=n$. 
Let $\mathcal A$ be an ample 
invertible sheaf on $W$ such that $|\mathcal A|$ is 
free. 
Then $\mathcal F\otimes \mathcal A^{\otimes n+1}$ 
{\em{(}}resp.~$\mathcal F\otimes \mathcal A^{\otimes n}${\em{)}} 
is generated by global sections. 
\end{lem}
\begin{proof}
This is a direct consequence 
of Lemmas \ref{p-lem5.7}, \ref{p-lem7.1}, and 
Castelnuovo--Mumford regularity. 
\end{proof}

Let us recall a vanishing theorem for dlt pairs. 

\begin{lem}\label{p-lem7.4}
Let $f\colon V\to W$ be a surjective morphism 
from a smooth 
projective variety $V$ onto a normal projective 
variety $W$. 
Let $\Delta$ be an effective $\mathbb R$-divisor 
on $V$ such that 
$(V, \Delta)$ is dlt and that every log canonical center of $(V, \Delta)$ 
is dominant onto $W$. 
Let $L$ be a Cartier divisor on $V$ such that 
$L-(K_V+\Delta)\sim _{\mathbb R} f^*H$ for some nef 
and big $\mathbb R$-divisor $H$ on $W$. 
Then $H^i(W, R^jf_*\mathcal O_V(L)\otimes \mathcal N)=0$ for $i>0$, 
$j\geq 0$, and every nef invertible sheaf $\mathcal N$ on $W$. 
\end{lem}

\begin{proof}[Sketch of Proof]
By Kodaira's lemma, we can write $H\sim _{\mathbb R}A+E$ such that 
$A$ is an ample $\mathbb R$-divisor on $W$ and 
$E$ is an effective $\mathbb R$-Cartier 
$\mathbb R$-divisor 
on $W$. 
Since every log canonical center of $(V, \Delta)$ is dominant onto $W$, 
$(V, \Delta+\varepsilon f^*E)$ is dlt for 
$0<\varepsilon \ll 1$. 
Let $N$ be a Cartier divisor on $W$ such that 
$\mathcal N\simeq \mathcal O_W(N)$. 
We note that  
$$
L+f^*N-(K_V+\Delta+\varepsilon f^*E)\sim _{\mathbb R} 
f^*(N+(1-\varepsilon )H+\varepsilon A)
$$
and that $N+(1-\varepsilon )H+\varepsilon A$ is ample for 
$0<\varepsilon \ll 1$. 
By \cite[Lemma 7.14]{fujino-weak}, 
$$H^i(W, R^jf_*\mathcal O_V(L)\otimes \mathcal N)=0$$ 
for $i>0$ and $j\geq 0$. 
\end{proof}

In Lemmas \ref{p-lem7.5} and \ref{p-lem7.6}, 
we treat mixed-$\widehat\omega$-big-sheaves on smooth 
projective curves. 

\begin{lem}\label{p-lem7.5}
Let $\mathcal G$ be a mixed-$\widehat\omega$-big-sheaf 
on a smooth projective curve $C$. 
Then $H^1(C, \mathcal G\otimes \mathcal N)=0$ holds for 
every nef invertible sheaf 
$\mathcal N$ on $C$.  
\end{lem}

\begin{proof}
We note that $\mathcal G$ is locally free since $C$ is a smooth curve. 
By definition, we have a mixed-$\omega$-big-sheaf 
$\mathcal F$ on $C$ and a 
generically isomorphic injection $\iota\colon  \mathcal F\hookrightarrow 
\mathcal G$. Note that the cokernel of $\iota$ is a skyscraper sheaf 
on $C$. 
By Lemma \ref{p-lem7.1}, 
$H^1(C, \mathcal F\otimes \mathcal N)=0$ holds. 
Therefore, we have $H^1(C, \mathcal G\otimes \mathcal N)=0$ 
by the surjection 
$H^1(C, \mathcal F\otimes \mathcal N)\to H^1(C, \mathcal G
\otimes \mathcal N)$. 
\end{proof}

\begin{lem}\label{p-lem7.6} 
Let $\mathcal E$ be a locally free sheaf on a smooth 
projective curve $C$ and let $P$ be a closed point of $C$. 
If $\mathcal E\otimes \mathcal O_C(-P)\otimes \mathcal N^{\otimes -1}$ is a 
mixed-$\widehat \omega$-big-sheaf 
on $C$ for some nef invertible sheaf $\mathcal N$ on $C$, then 
$\mathcal E$ is generated by global sections at $P$. 
\end{lem}

\begin{proof}
By Lemma \ref{p-lem7.5}, 
$H^1(C, \mathcal E\otimes \mathcal O_C(-P))=0$. 
This means that the natural restriction map 
$$
H^0(C, \mathcal E)\to \mathcal E\otimes \mathbb C(P)
$$ 
is surjective. Therefore, $\mathcal E$ is generated by 
global sections at $P$. 
\end{proof}

Let us discuss generically global generations of 
mixed-$\omega$-big-sheaves. 

\begin{lem}\label{p-lem7.7}
Let $\mathcal F$ be a mixed-$\omega$-sheaf on a normal 
projective variety $W$ with 
$\dim W=n$. 
Let $H$ be a big Cartier divisor on $W$ such that 
$|H|$ is free. 
Then $\mathcal F\otimes \mathcal O_W((n+1)H)$ is 
generically generated by global sections. 
\end{lem}
\begin{proof} 
If $W$ is a curve, then 
$H$ is ample. Therefore, 
the statement follows from Lemma \ref{p-lem7.3} when $n=1$. 
We will use induction on $n$. 
We may assume that $\mathcal F=
f_*\mathcal O_V(K_V+D)$ as in Definition \ref{p-def5.1}. 
\setcounter{step}{0}
\begin{step}\label{p-step7.7.1}
Let $\mu\colon \widetilde V\to V$ be a projective 
birational morphism 
from a smooth 
projective variety $\widetilde V$ such that 
$$
K_{\widetilde V}+\widetilde D=\mu^*(K_V+D)+E
$$ 
where $\widetilde D$ and $E$ are effective 
divisors and have no common irreducible components. 
Since 
$\mu_*\mathcal O_{\widetilde V}(K_{\widetilde V}
+\widetilde D)\simeq 
\mathcal O_V(K_V+D)$, we may replace 
$(V, D)$ and $f\colon V\to W$ with $(\widetilde V, \widetilde D)$ and 
$f\circ\mu\colon  \widetilde V\to W$, respectively. 
By taking $\mu\colon \widetilde V\to V$ suitably, 
we may assume that all the log canonical 
centers of $(V, D_{\mathrm{hor}})$ are 
dominant onto $W$, where $D_{\mathrm{hor}}$ is 
the horizontal part of $D$. 
Since $f_*\mathcal O_V(K_V+D_{\mathrm{hor}})
\hookrightarrow f_*\mathcal O_V(K_V+D)$ is a 
generically isomorphic injection, we may 
replace $D$ with $D_{\mathrm{hor}}$. 
\end{step}
\begin{step}\label{p-step7.7.2}
We will prove that 
$f_*\mathcal O_V(K_V+D)\otimes 
\mathcal O_V((n+1)H)$ is generically 
generated by global sections by induction on $n=\dim W$. 
We take a general member $W'$ of $|H|$. We put $f^{-1}(W')=V'$. 
Then we have a short exact sequence 
$$
0\to \mathcal O_V(K_V+D)\to \mathcal O_V(K_V+V'+D)\to 
\mathcal O_{V'}(K_{V'}+D|_{V'})\to 0 
$$ 
by adjunction. Since $W'$ is a general member of $|H|$, 
$$R^1f_*\mathcal O_V(K_V+D)\otimes 
\mathcal O_W(nH)\to R^1f_*\mathcal O_V(K_V+V'+D)\otimes 
\mathcal O_W(nH)$$ is injective. Hence we get a 
short exact sequence 
\begin{equation}\label{p-7.1}
\begin{split}
0&\to f_*\mathcal O_V(K_V+D)\otimes 
\mathcal O_W(nH)\to f_*\mathcal O_V(K_V+V'+D)\otimes 
\mathcal O_W(nH)\\&\to 
f_*\mathcal O_{V'}(K_{V'}+D|_{V'})\otimes 
\mathcal O_{W'}(nH|_{W'})\to 0. 
\end{split}
\end{equation}
By the vanishing theorem (see Lemma \ref{p-lem7.4}), 
we have 
\begin{equation}\label{p-7.2}
H^1(W, f_*\mathcal O_V(K_V+D)\otimes 
\mathcal O_W(nH))=0. 
\end{equation} 
Therefore, the restriction map 
$$
H^0(W, f_*\mathcal O_V(K_V+D)\otimes 
\mathcal O_W((n+1)H))\to 
H^0(W', f_*\mathcal O_{V'}(K_{V'}+D|_{V'})\otimes 
\mathcal O_{W'}(nH|_{W'}))
$$
is surjective by \eqref{p-7.1} and \eqref{p-7.2}. 
By induction on $n$, 
$f_*\mathcal O_{V'}(K_{V'}+D|_{V'})
\otimes \mathcal O_{W'}(nH|_{W'})$ is generically generated by global sections. 
This implies that 
so is $f_*\mathcal O_V(K_V+D)\otimes 
\mathcal O_W((n+1)H)$. 
\end{step} 
We obtain the desired statement. 
\end{proof}

Lemma \ref{p-lem7.8} is similar to Lemma \ref{p-lem7.7}. 

\begin{lem}\label{p-lem7.8}
Let $\mathcal F$ be a mixed-$\omega$-big-sheaf on a normal 
projective variety $W$ with 
$\dim W=n$. 
Let $H$ be a big Cartier divisor on $W$ such that 
$|H|$ is free. 
Then $\mathcal F\otimes \mathcal O_W(nH)$ 
is generically generated by global sections. 
\end{lem}

The proof of Lemma \ref{p-lem7.8} is essentially the same 
as that of Lemma \ref{p-lem7.7}. 

\begin{proof}[Sketch of Proof of Lemma \ref{p-lem7.8}]
If $n=0$, then the statement is trivial. 
If $n=1$, then it follows from Lemma \ref{p-lem7.6}. 
Therefore, we assume that $n\geq 2$. 
As in the proof of Lemma \ref{p-lem7.7}, 
we may assume that $\mathcal F=f_*\mathcal O_V(K_V+D+P)$ as in 
Definition \ref{p-def5.3}. 
Moreover, we may further assume that 
$D=D_{\mathrm{hor}}$ and that every log canonical 
center of $(V, D)$ is dominant onto $W$ (see 
Step \ref{p-step7.7.1} in the proof of Lemma \ref{p-lem7.7}). 
We take a general 
member $W'$ of $|H|$ and put $V'=f^{-1}(W')$. 
Then the natural restriction map 
\begin{equation*}
\begin{split}
&H^0(W, f_*\mathcal O_V(K_V+D+P)\otimes 
\mathcal O_W(nH))\\ &\to 
H^0(W', f_*\mathcal O_{V'}(K_{V'}+D|_{V'}+P|_{V'})\otimes 
\mathcal O_{W'}((n-1)H|_{W'}))
\end{split}
\end{equation*} 
is surjective as in 
Step \ref{p-step7.7.2} in the proof of 
Lemma \ref{p-lem7.7}. 
By induction on dimension, we see that 
$$
f_*\mathcal O_V(K_V+D+P)\otimes \mathcal O_W(nH)
$$ 
is generically generated by global sections.  
\end{proof}

We close this section with the following 
result, which is due to \cite{dutta-murayama}. 
We will use it in the proof of Theorem \ref{q-thm1.9}. 

\begin{lem}\label{p-lem7.9} 
Let $\mathcal F$ be a mixed-$\omega$-sheaf on a normal projective 
variety $W$ and let $H$ be a nef and big Cartier divisor on $W$. 
We put $\dim W=n$. 
Then $\mathcal F\otimes \mathcal O_W(lH)
$ is generically generated by global sections for $l\geq n^2+1$. 
\end{lem}
\begin{proof}
We may assume that 
$\mathcal F=f_*\mathcal O_V(K_V+D)$ as in Definition \ref{p-def5.1}. 
Then, by \cite[Theorem 1]{ekl} and 
\cite[Theorems C and 2.20]{dutta-murayama}, 
$\mathcal F\otimes \mathcal O_W(lH)
$ is generically generated by global sections for $l\geq n^2+1$. 
\end{proof}

For the details of Lemma \ref{p-lem7.9}, we recommend the 
reader to see \cite{dutta-murayama}. 

\section{A special case}\label{p-sec8} 

In this section, we freely use the standard notation and 
some basic results in the theory of minimal models 
(see, for example, \cite{fujino-fundamental}, \cite{fujino-foundations}, 
and \cite{fujino-iitaka}). 
We treat weakly semistable morphisms 
$f\colon X\to Y$ in the sense of Abramovich--Karu with 
the assumption that the geometric generic fiber $X_{\overline \eta}$ 
of $f\colon X\to Y$ has 
a good minimal model. 
In this case, we can prove some strong results with 
the aid of the theory of minimal models. We do not need the results of 
this section in the subsequent sections. Hence 
the reader can skip this section. 

\medskip 

Here we adopt the following definition of weakly semistable morphisms 
in the sense of Abramovich--Karu (see \cite{abramovich-karu}). 

\begin{defn}[Weakly semistable morphisms]\label{p-def8.1} 
Let $f\colon X\to Y$ be a projective surjective 
morphism between normal 
quasi-projective varieties with connected fibers. 
We say that $f\colon X\to Y$ is {\em{weakly semistable}} 
if 
\begin{itemize}
\item[(i)] the varieties $X$ and $Y$ admit toroidal structures 
$(U_X\subset X)$ and $(U_Y\subset Y)$ with 
$U_X=f^{-1}(U_Y)$, 
\item[(ii)] with this structure, the morphism $f$ is toroidal, 
\item[(iii)] the morphism $f$ is equidimensional, 
\item[(iv)] all the fibers of the morphism $f$ are reduced, 
and 
\item[(v)] $Y$ is smooth.  
\end{itemize}
\end{defn}

The following result is the main theorem of this section. 

\begin{thm}\label{p-thm8.2}
Let $f\colon X\to Y$ be a surjective 
morphism from a normal projective 
variety $X$ onto a smooth projective 
variety $Y$ with connected fibers. 
Assume that $f$ is weakly semistable in the sense of 
Abramovich--Karu and that the geometric generic 
fiber $X_{\overline \eta}$ of $f\colon  X\to Y$ has a good minimal model. 
Let $H$ be an ample Cartier divisor on $Y$. 
Let $k$ be a positive integer such that 
$k\geq 2$ and $f_*\omega^{\otimes k}_{X/Y}\ne 0$. 
Then 
$$
f_*\omega^{\otimes k}_{X/Y}\otimes \omega_Y\otimes \mathcal 
O_Y(H)
$$ 
is locally free and is a pure-$\omega$-big-sheaf on $Y$. 
More generally, we obtain that 
$$
\left(\bigotimes ^s f_*\omega^{\otimes k}_{X/Y}\right)\otimes \omega_Y 
\otimes \mathcal O_Y(H)
$$ 
is a pure-$\omega$-big-sheaf on $Y$ for every positive integer $s$. 
Therefore, if $A$ is an ample Cartier divisor on $Y$ such that 
$|A|$ is free, then 
$$
\left(\bigotimes ^s f_*\omega^{\otimes k}_{X/Y}\right)\otimes \omega_Y 
\otimes \mathcal O_Y(H+nA)
$$ 
is generated by global sections, where $n=\dim Y$. 
\end{thm}

\begin{proof}
As mentioned above, we will freely use some basic results in the theory 
of minimal models. 
We note that $X$ has only rational Gorenstein singularities by 
\cite[Lemma 6.1]{abramovich-karu} since 
$f\colon X\to Y$ is weakly semistable by assumption. 
Hence $X$ has only canonical Gorenstein singularities. 
\setcounter{step}{0}
\begin{step}\label{p-step8.2-1}
By the proof of \cite[Theorem 1.6]{fujino-direct} 
(see also \cite{fujino-direct-corri}), 
we have already known that 
$f_*\omega^{\otimes m}_{X/Y}$ is a nef locally free sheaf on 
$Y$ for every $m\geq 1$. 
\end{step}
\begin{step}\label{p-step8.2-2}
We consider a relative good minimal model $f'\colon X'\to Y$ 
of $f\colon X\to Y$ (see \cite[Theorem 3.3]{fujino-direct}). 
$$
\xymatrix{
X \ar[dr]_-f\ar@{-->}^-\phi[rr]&& X'\ar[dl]^-{f'} \\ 
& Y& 
}
$$ 
Since 
$$
f_*\omega^{\otimes m}_{X/Y}\simeq f'_*\mathcal O_{X'}(mK_{X'/Y})
$$ 
holds for every $m\geq 1$, it is sufficient to 
prove that 
$$
f'_*\mathcal O_{X'}(K_{X'}+(k-1)K_{X'/Y}+f'^*H)
$$ 
is a pure-$\omega$-big-sheaf on $Y$. 
\end{step}
\begin{step}\label{p-step8.2-3}
In this step, we will prove: 
\begin{claim}
$K_{X'/Y}$ is nef and $f'$-semi-ample. 
\end{claim}
\begin{proof}[Proof of Claim]
Since $f'\colon  X'\to Y$ is a relative good minimal model 
of $f\colon X\to Y$, $K_{X'/Y}$ is $f'$-semi-ample. 
Therefore, 
$$
f'^*f'_*\mathcal O_{X'}(lK_{X'/Y})\to \mathcal O_{X'}(lK_{X'/Y})
$$ 
is surjective for a sufficiently large and divisible 
positive integer $l$. 
Since $f'_*\mathcal O_{X'}(lK_{X'/Y})\simeq 
f_*\omega^{\otimes l}_{X/Y}$ is a nef locally 
free sheaf, $K_{X'/Y}$ is nef by the above 
surjection. 
\end{proof}
\end{step}
\begin{step}\label{p-step8.2-4}
Since $K_{X'/Y}$ is nef and $f'$-semi-ample, 
$(k-1)K_{X'/Y}+af'^*H$ is semi-ample 
for every positive 
rational number $a$ by Lemma \ref{p-lem3.5}. 
Since $X'$ is a relative minimal model of $f\colon X\to Y$ and 
$X$ has only canonical singularities, 
$X'$ also has only canonical singularities. 
We take a birational morphism 
$\rho\colon \widetilde X\to X'$ from a smooth 
projective variety $\widetilde X$ such that 
the exceptional locus $\Exc(\rho)$ of $\rho$ 
is a simple normal crossing 
divisor on $\widetilde X$. 
Since $X'$ has only canonical singularities, 
we see that 
$$
\rho_*\mathcal O_{\widetilde X}(K_{\widetilde X}+\lceil 
(k-1)\rho^*K_{X'/Y}+\rho^*f'^*H\rceil)
\simeq \mathcal O_{X'}(K_{X'}+(k-1)K_{X'/Y}+f'^*H)
$$ 
holds. By Lemma \ref{p-lem5.9}, 
$$
(f'\circ\rho)_*\mathcal O_{\widetilde X}(K_{\widetilde X}+\lceil 
(k-1)\rho^*K_{X'/Y}+\rho^*f'^*H\rceil)
\simeq f'_*\mathcal O_{X'}(K_{X'}+(k-1)K_{X'/Y}+f'^*H)
$$ is a pure-$\omega$-big-sheaf on $Y$ since 
$$
(k-1)\rho^*K_{X'/Y}+\rho^*f'^*H-\frac{1}{2}\rho^*f'^*H
$$ 
is semi-ample. 
This means that 
$$
f_*\omega^{\otimes k}_{X/Y}\otimes \omega_Y\otimes \mathcal 
O_Y(H)
$$ 
is locally free and a pure-$\omega$-big-sheaf on $Y$. 
\end{step}
\begin{step}\label{p-step8.2-5} 
In this step, we will briefly explain how to prove 
that $$
\left(\bigotimes ^s f_*\omega^{\otimes k}_{X/Y}\right)\otimes \omega_Y 
\otimes \mathcal O_Y(H)
$$ 
is a pure-$\omega$-big-sheaf on $Y$ for every positive integer $s$. 

Let $$
X^s=\underbrace{X\times _Y X \times_Y \cdots \times _Y X}_{s}
$$ 
be the $s$-fold fiber product of $f\colon X\to Y$ and let 
$f^s\colon X^s\to Y$ be the induced morphism. 
Then we can check that 
$X^s$ has only rational Gorenstein singularities since $f\colon X\to Y$ 
is weakly semistable (see \cite[Section 5]{fujino-direct} and 
Step 4 in the proof of \cite[Theorem 4.3.1]{fujino-iitaka}). 
Therefore, $X^s$ has only canonical Gorenstein singularities.  
Let $F$ be a general fiber of $f\colon X\to Y$. 
Then $F$ has a good minimal model by assumption and 
\cite[Theorem 3.3]{fujino-direct}. This implies that 
$$
F^s=\underbrace{F\times F\times \cdots \times F}_{s}
$$ 
also has a good minimal model. 
We note that $F^s$ is a general fiber of $f^s\colon X^s\to Y$. 
Hence $f^s\colon  X^s\to Y$ has a relative good minimal model over $Y$ 
by \cite[Theorem 3.3]{fujino-direct}. 
By the flat base change theorem and the projection formula, 
we can check that 
\begin{equation}\label{p-eq8.1}
\left(\bigotimes ^s f_*\omega^{\otimes k}_{X/Y}\right)\simeq 
f^s_*\omega^{\otimes k}_{X^s/Y}
\end{equation}
holds (see \cite[Section 5]{fujino-direct} and 
Step 4 in the proof of \cite[Theorem 4.3.1]{fujino-iitaka}). 
In particular, $f^s_*\omega^{\otimes k}_{X^s/Y}$ is a nef 
locally free sheaf on $Y$. 
By applying the arguments in 
Steps \ref{p-step8.2-2}, \ref{p-step8.2-3}, and 
\ref{p-step8.2-4} to $f^s\colon X^s\to Y$, 
we obtain that 
\begin{equation}\label{p-eq8.2}
f^s_*\omega^{\otimes k}_{X^s/Y}\otimes \omega_Y \otimes 
\mathcal O_Y(H)
\end{equation} 
is a pure-$\omega$-big-sheaf on $Y$. 
Therefore, by \eqref{p-eq8.1} and \eqref{p-eq8.2}, we have that 
$$
\left(\bigotimes ^s f_*\omega^{\otimes k}_{X/Y}\right)\otimes \omega_Y 
\otimes \mathcal O_Y(H)
$$ 
is a pure-$\omega$-big-sheaf on $Y$ 
for every positive integer $s$. 
\end{step}
By Lemma \ref{p-lem7.3}, 
$$
\left(\bigotimes ^s f_*\omega^{\otimes k}_{X/Y}\right)\otimes \omega_Y 
\otimes \mathcal O_Y(H+nA)
$$ 
is generated by global sections. 
\end{proof}

\begin{rem}\label{p-rem8.3} 
In the proof of \cite[Theorem 3.3]{fujino-direct}, 
which was used in the proof of Theorem \ref{p-thm8.2}, 
we need Theorem \ref{q-thm1.13} (see \cite[Theorem 3.4]{fujino-direct}). 
We note that we will not use Theorem \ref{p-thm8.2} in the subsequent 
sections. 
Hence there is no problem. 
\end{rem}

We note that the global generation of 
$$
\left(\bigotimes ^s f_*\omega_{X/Y}\right)\otimes \omega_Y 
\otimes \mathcal O_Y(H+nA)
$$ was already treated in \cite[Theorem 3.6]{kollar}. 
In some sense, Theorem \ref{p-thm8.2} generalizes 
\cite[Theorem 1.8]{su-yang}. 

Theorem \ref{p-thm8.2} predicts that $f_*\omega^{\otimes k}
_{X/Y}\otimes \omega_Y 
\otimes \mathcal O_Y(H)$ has good properties. 
Of course, we strongly hope to prove Theorem 
\ref{p-thm8.2} without using the assumption that 
$X_{\overline \eta}$ has a good minimal model. 

\section{Fundamental theorem}\label{p-sec9}

This section is the main part of this paper. 
The main result of this section is Theorem \ref{p-thm9.3}, 
which we call a fundamental theorem of the theory of 
mixed-$\omega$-sheaves.  

\medskip 

Let us start with the following lemma. 

\begin{lem}[{\cite[Chapter V, 3.34.~Lemma]{nakayama}}]\label{p-lem9.1}
Let $f\colon X\to Y$ be a surjective morphism from a normal projective 
variety $X$ onto a smooth projective variety $Y$. 
Let $L$ be a Cartier divisor on $X$ and 
let $\Delta$ be an effective $\mathbb R$-divisor 
on $X$ such that 
$K_X+\Delta$ is $\mathbb R$-Cartier. 
Let $k$ be a positive integer with $k\geq 2$. 
We assume the following conditions:
\begin{itemize}
\item[(i)] $(X, \Delta)$ is log canonical over a nonempty 
Zariski open set of $Y$, and 
\item[(ii)] $L-k(K_{X/Y}+\Delta)$ is nef and 
$f$-semi-ample. 
\end{itemize}
Let $H$ be an ample divisor on $Y$. 
We assume that $f_*\mathcal O_Y(L)\ne 0$. 
We take a positive integer $l$ such that 
$$
\mathcal O_Y(lH)\otimes f_*\mathcal O_X(L)
$$ 
is big. Then 
$$
\mathcal O_Y(K_Y+(l-\lfloor l/k\rfloor)H)\otimes 
f_*\mathcal O_X(L)
$$ 
is a mixed-$\widehat \omega$-big-sheaf 
on $Y$. 
Hence we obtain that 
$$
\mathcal O_Y(K_Y+(k-1)H)\otimes f_*\mathcal O_X(L)
$$ 
is always a mixed-$\widehat \omega$-big-sheaf 
on $Y$. 
\end{lem}

We include all the details although 
Lemma \ref{p-lem9.1} is essentially the same as 
\cite[Chapter V, 3.34.~Lemma]{nakayama}. 

\begin{proof}[Proof of Lemma \ref{p-lem9.1}]
We divide the proof into several small steps. 
\setcounter{step}{0} 
\begin{step}[Resolution of singularities]\label{p-step9.1.1}
Let $\mu\colon \widetilde X\to X$ be a projective birational morphism 
from a smooth projective variety $\widetilde X$ such that 
$K_{\widetilde X}+\widetilde \Delta=\mu^*(K_X+\Delta)$ and 
that $\Supp \widetilde \Delta$ is a simple normal crossing 
divisor on $\widetilde X$. 
We put $E=\lceil -(\widetilde \Delta^{<0})\rceil$. 
Then $E$ is an effective $\mu$-exceptional divisor on $\widetilde X$, 
$\widetilde \Delta+E$ is effective, 
and $(\widetilde X, \widetilde \Delta+E)$ is log canonical over 
a nonempty Zariski open set of $Y$ by construction. 
We note that 
$$
\mu^*L+kE-k(K_{\widetilde X/Y}+\widetilde \Delta+E)
=\mu^*(L-k(K_{X/Y}+\Delta))
$$ 
and that $\mu_*\mathcal O_{\widetilde X}(\mu^*L+kE)\simeq 
\mathcal O_X(L)$. 
Therefore, by replacing $f\colon X\to Y$, $L$, and $\Delta$ with 
$f\circ \mu\colon  \widetilde X\to Y$, $\mu^*L+kE$, and 
$\widetilde \Delta+E$ respectively, 
we may assume that 
$X$ is smooth and $\Supp \Delta$ is a simple normal crossing divisor 
on $X$. 
\end{step}
\begin{step}\label{p-step9.1.2} 
We note that 
by Lemma \ref{p-lem4.6} we can 
always take a positive integer $l$ such 
that $\mathcal O_Y(lH)\otimes 
f_*\mathcal O_X(L)$ is big since 
$H$ is an ample divisor on $Y$. 
Since $\mathcal O_Y(lH)\otimes f_*\mathcal O_X(L)$ is big, 
we can take a positive integer $a$ such that 
$$\widehat S^a(\mathcal O_Y(lH)\otimes f_*\mathcal O_X(L))
\otimes \mathcal O_Y(-H)=\mathcal O_Y((al-1)H)\otimes 
\widehat S^a(f_*\mathcal O_X(L))$$ is generically generated by 
global sections by Lemma \ref{p-lem4.4}. 
\end{step}
\begin{step}\label{p-step9.1.3}
We take an effective $f$-exceptional divisor $E$ on $X$ such that 
$$
\left(f_*\mathcal O_X(bL)\right)^{**}\simeq f_*\mathcal O_X(b(L+E))
$$
holds for every $1\leq b\leq a$. 
For a related result, see \cite[Theorem 1.2]{fujino-nakayama}. 
By taking a 
resolution of singularities as in Step \ref{p-step9.1.1}, we may assume that 
$\Supp (\Delta+E)$ is a simple normal crossing divisor on $X$. 
Since 
$$
(L+E)-k(K_{X/Y}+\Delta+(1/k)E)=L-k(K_{X/Y}+\Delta), 
$$ 
we may replace $L$ and $\Delta$ with 
$L+E$ and $\Delta+(1/k)E$, respectively. 
This is because 
$$
\mathcal O_Y(K_Y+(l-\lfloor l/k\rfloor)H)\otimes 
f_*\mathcal O_X(L)
$$ 
is a mixed-$\widehat\omega$-big-sheaf on $Y$ if and only if 
so is 
$$\mathcal O_Y(K_Y+(l-\lfloor l/k\rfloor)H)\otimes 
\left(f_*\mathcal O_X(L)\right)^{**}. 
$$ 
Hence we may assume that 
$f_*\mathcal O_X(bL)$ is reflexive for every $1\leq b\leq a$. 
\end{step}
\begin{step}\label{p-step9.1.4}
By taking a suitable birational modification of $X$ again (see Step 
\ref{p-step9.1.1}), 
we may further assume that 
the image of the natural map 
$$
f^*f_*\mathcal O_X(L)\to \mathcal O_X(L)
$$ 
is invertible and can be written 
as $\mathcal O_X(L-B)$ such that 
$\Supp (\Delta+B)$ is a simple normal crossing divisor 
on $X$. 
By the definition of $B$, 
we have $f_*\mathcal O_X(L-B)=f_*\mathcal O_X(L)$. 
\end{step}
\begin{step}\label{p-step9.1.5}
We note that we can take an effective $f$-exceptional divisor $E$ on $X$ such that 
the map 
$f^*f_*\mathcal O_X(L)\to \mathcal O_X(L-B)$ induces 
$$
f^*\widehat S^a(f_*\mathcal O_X(L))\to \mathcal O_X(a(L-B)+E). 
$$
Then we have the following map 
\begin{equation}\label{p-9.1}
H^0(Y, \mathcal O_Y((al-1)H)\otimes 
\widehat S^a(f_*\mathcal O_X(L)))\otimes \mathcal O_X\to 
\mathcal O_X(a(L-B)+E+(al-1)f^*H). 
\end{equation}
By taking a suitable birational modification of $X$ again 
(see Step \ref{p-step9.1.1}), we may assume that 
the image of \eqref{p-9.1} 
is $$
\mathcal O_X(a(L-B)+E-F+(al-1)f^*H)
$$ 
for some effective $f$-vertical divisor $F$ on $X$. 
We may further assume that 
$\Supp (\Delta+B+E+F)$ is a simple normal crossing divisor 
on $X$. 
We put 
$$
N:=a(L-B)+E-F+(al-1)f^*H. 
$$ 
Then $|N|$ is free by \eqref{p-9.1} and the definition of $F$. 
\end{step}
\begin{step}\label{p-step9.1.6} 
We take a positive number $\varepsilon$. 
Then, by Lemma \ref{p-lem3.5}, we obtain 
that $L-k(K_{X/Y}+\Delta)+\varepsilon f^*H$ is 
semi-ample because 
$L-k(K_{X/Y}+\Delta)$ is nef and $f$-semi-ample by assumption. 
We put 
$$
M:=L-(K_{X/Y}+\Delta)-\frac{k-1}{k}B+\frac{k-1}{ak}(E-F)
+\left(l-\left\lfloor \frac{l}{k}\right\rfloor\right)f^*H. 
$$ 
We note 
that 
$$
\frac{(al-1)(k-1)}{ak}<\left\lceil\frac{l(k-1)}{k}\right\rceil
=l-\left\lfloor\frac{l}{k}\right\rfloor. 
$$ 
Then 
$$
M-\frac{k-1}{ak}N-\frac{1}{k}(L-k(K_{X/Y}+\Delta)+\varepsilon f^*H)=\alpha f^*H
$$ 
for some $\alpha>0$ if $\varepsilon$ is sufficiently small. 
Thus $M$ and $M-\alpha f^*H$ are 
semi-ample. Without loss of generality, we may assume that 
$\varepsilon$ and $\alpha$ are rational numbers since 
$$
\left(l-\left\lfloor \frac{l}{k}\right\rfloor\right)-\frac{(al-1)(k-1)}{k} 
-\frac{\varepsilon}{k}=\alpha. 
$$ 
\end{step}
\begin{step}\label{p-step9.1.7}
We consider 
$$
\left\lfloor \frac{k-1}{k}B+\Delta\right\rfloor. 
$$
We put 
$$
B_0=\max\left\{ T \left| \text{$T$ 
is a Weil divisor 
with $0\leq T\leq B$ and $T\leq 
\left\lfloor \frac{k-1}{k}B+
\Delta\right\rfloor$}\right.\right\}.
$$ 
We write 
$$
\left\lfloor \frac{k-1}{k}B+\Delta\right\rfloor-B_0=\Delta_1+\Delta_2
$$ 
where $\Delta_1$ is the horizontal part and $\Delta_2$ is the 
vertical part. 
We note the following easy Claim. 
\begin{claim}\label{p-claim9.1} 
Let $r$ be a real number with $0\leq r\leq 1$, 
let $k$ be a positive integer with $k\geq 2$, and let 
$b$ be a nonnegative integer. 
Then 
$$
\left\lfloor \frac{k-1}{k}b+r\right\rfloor 
-\min \left\{ b, \left\lfloor \frac{k-1}{k}b+r\right\rfloor \right\}
=\begin{cases}
1 & \text{if $r=1$ and $b=0$}, \\ 
0 &\text{otherwise}. 
\end{cases}
$$
\end{claim}
\begin{proof}[Proof of Claim] 
We assume that 
$$
b\leq \left\lfloor \frac{k-1}{k}b+r\right\rfloor-1
$$ 
holds. Then 
$$
b\leq\frac{k-1}{k}b+r-1. 
$$ 
Hence we obtain $0\leq \frac{b}{k}\leq r-1\leq 0$. 
Thus we get $r=1$ and $b=0$. By this observation, 
we obtain the desired result. 
\end{proof}
By assumption (i) and the Claim, 
$\Delta_1=\Delta^{=1}_1$. 
By construction and the Claim, we see that 
$\Delta_1\subset \Supp \Delta^{=1}$ and 
that $\Delta_1$ and $\Supp \{M\}$ have no common irreducible 
components. 
\end{step}
\begin{step}\label{p-step9.1.8} 
We have the following generically isomorphic injections: 
\begin{equation*}
\begin{split}
f_*\mathcal O_X(K_X+\Delta_1+\lceil M\rceil)&\hookrightarrow 
\omega_Y((l-\lfloor l/k\rfloor)H)\otimes 
\left( f_*\mathcal O_X
\left(L-\left\lfloor\frac{k-1}{k}B+\Delta\right\rfloor+\Delta_1\right)\right)^{**}
\\& 
=\omega_Y((l-\lfloor l/k\rfloor)H)\otimes 
\left( f_*\mathcal O_X\left(L-B_0-\Delta_2\right)\right)^{**}
\\ &\hookrightarrow 
\omega_Y((l-\lfloor l/k\rfloor)H)\otimes 
f_*\mathcal O_X(L). 
\end{split}
\end{equation*}
We note that 
$$
f_*\mathcal O_X(L)=f_*\mathcal O_X(L-B)\subset 
f_*\mathcal O_X(L-B_0)\subset f_*\mathcal O_X(L). 
$$
This implies that 
$$
\mathcal O_Y(K_Y+(l-\lfloor l/k\rfloor)H)\otimes 
f_*\mathcal O_X(L)
$$ 
is a mixed-$\widehat \omega$-big-sheaf 
on $Y$ because $f_*\mathcal O_X(K_X+\Delta_1+\lceil M\rceil)$ 
is a mixed-$\omega$-big-sheaf by Lemma \ref{p-lem5.9}. 
\end{step}
\begin{step}\label{p-step9.1.9}
Let $l_0$ be the minimum positive integer such that 
$$
\mathcal O_Y(K_Y+l_0H)\otimes f_*\mathcal O_X(L)
$$ 
is a mixed-$\widehat\omega$-big-sheaf on $Y$. 
By Theorem \ref{p-thm6.3}, 
$$
\mathcal O_Y(l_0H)\otimes f_*\mathcal O_X(L)
$$ 
is big. By the result obtained above, 
$$
\mathcal O_Y(K_Y+(l_0-\lfloor l_0/k\rfloor)H)\otimes 
f_*\mathcal O_X(L)
$$ 
is a mixed-$\widehat\omega$-big-sheaf on $Y$. 
This implies that $l_0-\lfloor l_0/k\rfloor\geq l_0$. 
Thus we get $l_0\leq k-1$. 
Hence we have 
$$
\mathcal O_Y(K_Y+(k-1)H)\otimes f_*\mathcal O_X(L)
$$ 
is a mixed-$\widehat\omega$-big-sheaf on $Y$. 
\end{step}
Thus we get the desired statements. 
\end{proof}

\begin{rem}\label{p-rem9.2} In 
Lemma \ref{p-lem9.1}, 
we further assume that 
$(X, \Delta)$ is klt over a nonempty Zariski open set of $Y$. 
Then we can easily see that $\Delta_1=0$ 
in Step \ref{p-step9.1.7} in the proof of 
Lemma \ref{p-lem9.1}. 
Therefore, we obtain that 
$$
\mathcal O_Y(K_Y+(l-\lfloor l/k\rfloor)H)\otimes 
f_*\mathcal O_X(L) 
$$ 
and 
$$
\mathcal O_Y(K_Y+(k-1)H)\otimes 
f_*\mathcal O_X(L)
$$ 
are pure-$\widehat\omega$-big-sheaves on $Y$. 
\end{rem}

Theorem \ref{p-thm9.3} is the most important result 
in the theory of mixed-$\omega$-sheaves. 
So we call it a fundamental theorem of the theory 
of mixed-$\omega$-sheaves. 

\begin{thm}[{\cite[Chapter V, 3.35.~Theorem]{nakayama}}]\label{p-thm9.3} 
Let $f\colon X\to Y$ be a surjective morphism 
from a normal projective variety $X$ onto a smooth projective 
variety $Y$. 
Let $L$ be a Cartier divisor 
on $X$ and let $\Delta$ be an effective $\mathbb R$-divisor 
on $X$ such that 
$K_X+\Delta$ is $\mathbb R$-Cartier. 
Let $D$ be an $\mathbb R$-divisor on $Y$. 
Let $k$ be a positive integer with $k\geq 2$. 
Assume the following conditions: 
\begin{itemize}
\item[(i)] $(X, \Delta)$ is log canonical over a nonempty Zariski 
open set of $Y$, 
and 
\item[(ii)] $L+f^*D-k(K_{X/Y}+\Delta)-f^*A$ is 
semi-ample for some big $\mathbb R$-divisor 
$A$ on $Y$. 
\end{itemize}
If $f_*\mathcal O_Y(L)\ne 0$, then 
$$
\mathcal O_Y(K_Y+\lceil D\rceil)\otimes f_*\mathcal O_X(L)
$$ 
is a mixed-$\widehat \omega$-big-sheaf 
on $Y$. 
\end{thm}
\begin{proof}We divide the proof into several small steps. 
\setcounter{step}{0}
\begin{step}[Reductions]\label{p-step9.3.1}
By taking a resolution as in Step \ref{p-step9.1.1} in the proof of 
Lemma \ref{p-lem9.1}, we may assume that $X$ is a smooth projective 
variety and that $\Supp \Delta$ is a simple normal crossing divisor on 
$X$. 
We note that 
$$
L+f^*\lceil D\rceil -k\left(K_{X/Y}+\Delta+\frac{1}{k}
f^*\{-D\}\right)-f^*A=L+f^*D-k(K_{X/Y}+\Delta)-f^*A. 
$$ 
Therefore, by replacing $L$ and $\Delta$ with 
$L+f^*\lceil D\rceil$ and $\Delta+\frac{1}{k}f^*\{-D\}$, respectively, 
we may assume that $D=0$. 
By Kodaira's lemma, we have $A\sim_{\mathbb R} A_1+A_2$ such that 
$A_1$ is an ample $\mathbb R$-divisor and $A_2$ is an 
effective $\mathbb R$-divisor. 
By replacing $A$ and $\Delta$ with $A_1$ and 
$\Delta+\frac{1}{k}f^*A_2$ respectively, 
we may further assume that 
$A$ is an ample $\mathbb R$-divisor on $Y$. 
We take an ample Cartier divisor $H$ on $Y$ and a 
positive integer $m$ such that 
$A-\frac{k-1}{m}H$ is ample. 
Then 
$$
L-k(K_{X/Y}+\Delta)-\frac{k-1}{m}f^*H
$$ 
is semi-ample. 
Therefore, we may replace $A$ with 
$\frac{k-1}{m}H$. By taking a resolution 
as in Step \ref{p-step9.1.1} in the proof of Lemma 
\ref{p-lem9.1} again, we may assume that $X$ is a smooth 
projective variety and that $\Supp \Delta$ is 
a simple normal crossing 
divisor on $X$. 
By Lemma \ref{p-lem3.2}, 
we may further assume that $\Delta$ is a 
$\mathbb Q$-divisor. 
We take an effective $f$-exceptional 
divisor  $E$ and replace $L$ and $\Delta$ with 
$L+E$ and $\Delta+(1/k)E$ respectively. 
Then we may assume that 
$f_*\mathcal O_X(L)$ is reflexive. 
By taking a birational modification of $X$,  
we may assume that the image 
of 
$$
f^*f_*\mathcal O_X(L)\to \mathcal O_X(L)
$$ 
is $\mathcal O_X(L-B)$ for some effective 
divisor $B$ such that 
$\Supp (\Delta+B)$ is a simple normal crossing 
divisor on $X$. 
Let $S$ denote the union of all $f$-exceptional 
divisors on $X$. 
We may assume that 
$\Supp (\Delta+B+S)$ is a simple normal 
crossing divisor on $X$ by taking a suitable birational 
modification of $X$ again (see Step \ref{p-step9.1.1} 
in the proof of 
Lemma \ref{p-lem9.1}). 
\end{step}
\begin{step}\label{p-step9.3.2}
By Lemma \ref{p-lem3.3} and Remark \ref{p-rem3.4}, 
we take a finite flat Galois cover $\tau\colon Y'\to Y$ from a smooth 
projective variety $Y'$ and get the following 
commutative diagram 
$$
\xymatrix{
X' \ar[d]_-{f'}\ar[r]^-\rho& X\ar[d]^-f \\
Y' \ar[r]_-\tau&Y
}
$$
such that 
$X'=X\times _YY'$ is a smooth projective 
variety, 
$\tau^*H\sim mH'$ for some ample Cartier divisor 
$H'$ on $Y'$, and $\rho^*\omega^{\otimes n}_{X/Y}=\omega^{\otimes n}_{X'/Y'}$ for every integer $n$. 
Let $G$ denote the Galois group of $\tau\colon Y'\to Y$. 
By construction (see the proof of Lemma \ref{p-lem3.3}), 
we may assume that $H'$ is $G$-invariant. 
We put 
$L'=\rho^*L$, $B'=\rho^*B$, 
$\Delta'=\rho^*\Delta$, 
and $S'=\rho^*S$. 
Without loss of generality, 
we may assume that 
$\Supp (\Delta'+B'+S')$ 
is a simple normal crossing divisor on $X'$ and that 
$\rho^*(K_{X/Y}+\Delta)=K_{X'/Y'}+\Delta'$ holds 
(see Remark \ref{p-rem3.4}). 
We note that  
$$
L'-(k-1)f'^*H'-k(K_{X'/Y'}+\Delta')\sim 
_{\mathbb Q} \rho^*\left(L-k(K_{X/Y}+\Delta)-\frac{k-1}{m}
f^*H\right) 
$$ by construction. 
This implies that 
$\left(L'-(k-1)f'^*H'\right)-k(K_{X'/Y'}+\Delta')$ is semi-ample. 
\end{step}
\begin{step}\label{p-step9.3.3}
We apply Lemma \ref{p-lem9.1} to 
$$\left(L'-(k-1)f'^*H'\right)-k(K_{X'/Y'}+\Delta').$$ 
Then 
we obtain that $\mathcal O_{Y'}(K_{Y'})\otimes 
f'_*\mathcal O_{X'}(L')$ 
is a mixed-$\widehat\omega$-big-sheaf 
on $Y'$. 
Therefore, 
$f'_*\mathcal O_{X'}(L')$ is a big sheaf on $Y'$ by Theorem \ref{p-thm6.3}. 
Thus we can take a positive integer $a$ such that 
$\widehat S^a(f'_*\mathcal O_{X'}(L'))$ is generically 
generated by global sections (see Lemma \ref{p-lem4.4}). 
Then we take an effective $G$-invariant $f'$-exceptional 
divisor $E'$ on $X'$ such that 
$$
\left(f'_*\mathcal O_{X'}(bL')\right)^{**}\simeq 
f'_*\mathcal O_{X'}(b(L'+E'))
$$ 
holds for every $1\leq b\leq a$. 
By replacing $L'$, $\Delta'$, and $B'$ with $L'+E'$, $\Delta'+(1/k)E'$, 
and $B'+E'$ respectively, 
we may assume that 
$f'_*\mathcal O_{X'}(bL')$ is reflexive for every $1\leq b\leq a$. 
\end{step}
\begin{step}\label{p-step9.3.4}
We can take an effective $G$-invariant 
$f'$-exceptional divisor $E'$ on $X'$ such that 
the surjective map 
$$
f'^*f'_*\mathcal O_{X'}(L')\to \mathcal O_{X'}(L'-B')
$$ 
induces 
$$
f'^*\widehat S^a(f'_*\mathcal O_{X'}(L'))
\to \mathcal O_{X'}(a(L'-B')+E'). 
$$ 
Then we have the following map 
\begin{equation}\label{p-9.2} 
H^0(Y', \widehat S^a(f'_*\mathcal O_{X'}(L')))
\otimes \mathcal O_{X'}\to \mathcal O_{X'}
(a(L'-B')+E'). 
\end{equation}
By taking an equivariant resolution 
of singularities of $X'$, we may assume that 
the image of \eqref{p-9.2} is 
$$
\mathcal O_{X'}(a(L'-B')+E'-F')
$$ 
for some effective $G$-invariant 
$f'$-vertical divisor $F'$ on $X'$. Of course, 
we may assume that $\Supp(\Delta'+B'+E'+F')$ 
is a simple normal crossing divisor on $X'$. 
We put 
$$
N':=a(L'-B')+E'-F'. 
$$  
Then $|N'|$ is free by \eqref{p-9.2} and the definition of $F'$. 
We put 
$$
M':=L'-(K_{X'/Y'}+\Delta')-\frac{k-1}{k}B'+\frac{k-1}{ak}
(E'-F'). 
$$
Then $$
M'-\frac{k-1}{ak}N'-\frac{1}{k}(L'-k(K_{X'/Y'}+\Delta')-
(k-1)f'^*H')=\frac{k-1}{k}f'^*H'. 
$$ 
In particular, $M'$ and $M'-\frac{k-1}{k}f'^*H'$ are semi-ample. 
\end{step}
\begin{step}\label{p-step9.3.5}
We put 
$$
\left\lfloor \frac{k-1}{k}B'+\Delta'\right\rfloor=B'_0+\Delta'_1+\Delta'_2
$$ 
as in Step \ref{p-step9.1.7} in the proof of 
Lemma \ref{p-lem9.1}. 
Then $\Delta'_1$ is a $G$-invariant 
$f'$-horizontal simple normal crossing divisor on $X'$. 
As before, $\Supp \{M'\}$ and $\Delta'_1$ have no common irreducible components. 
Thus, by Lemma \ref{p-lem5.9}, 
$f'_*\mathcal O_{X'}(K_{X'}+\Delta'_1+\lceil M'\rceil)$ 
is a mixed-$\omega$-big-sheaf on $Y'$. 
Note that the Galois group $G$ acts on 
$f'_*\mathcal O_{X'}(K_{X'}+\Delta'_1+\lceil M'\rceil)$. 
\end{step}
\begin{step}\label{p-step9.3.6} 
Therefore, we get the following 
generically isomorphic $G$-equivariant embedding: 
\begin{equation}\label{p-9.3}
f'_*\mathcal O_{X'}(K_{X'}+\Delta'_1+\lceil M'\rceil)
\hookrightarrow \mathcal O_{Y'}(K_{Y'})
\otimes f'_*\mathcal O_{X'}(L') 
\end{equation} 
as in Step \ref{p-step9.1.8} in the proof of 
Lemma \ref{p-lem9.1}. 
We note that $f'_*\mathcal O_{X'}(L')\simeq 
\tau^*f_*\mathcal O_X(L)$ by 
the flat base change theorem. 
We take $\tau_*$ of \eqref{p-9.3} and 
then take the $G$-invariant parts. Thus, 
we get a mixed-$\omega$-big-sheaf 
$$
\mathcal F:=\left(\tau_*f'_*\mathcal O_{X'}(K_{X'}
+\Delta'_1+\lceil M'\rceil)\right)^G
$$ 
on $Y$ and a generically isomorphic injection 
$$
\mathcal F\hookrightarrow \mathcal O_Y(K_Y)\otimes 
f_*\mathcal O_X(L). 
$$ 
This means that $\mathcal O_Y(K_Y)\otimes 
f_*\mathcal O_X(L)$ is a mixed-$\widehat\omega$-big-sheaf on $Y$. 
\end{step} 
Hence we obtain that 
$\mathcal O_Y(K_Y+\lceil D\rceil)\otimes 
f_*\mathcal O_X(L)$ is a mixed-$\widehat\omega$-big-sheaf on $Y$. 
\end{proof}

\begin{rem}\label{p-rem9.4}
As in Remark \ref{p-rem9.2}, 
we further assume that 
$(X, \Delta)$ is klt over a nonempty Zariski open set of 
$Y$ in Theorem \ref{p-thm9.3}. 
Then we see that $\Delta'_1=0$ in Step 
\ref{p-step9.3.5} in the proof of Theorem \ref{p-thm9.3}. 
Hence we obtain that 
$$
\mathcal O_Y(K_Y+\lceil D\rceil)\otimes f_*\mathcal O_X(L) 
$$ 
is a pure-$\widehat\omega$-big-sheaf on $Y$. 
\end{rem}

As a corollary of Theorem \ref{p-thm9.3}, we have: 

\begin{cor}[{\cite[Chapter V.~3.37.~Corollary]{nakayama}}]
\label{p-cor9.5}
Let $f\colon X\to Y$ be a surjective morphism from a normal projective 
variety $X$ onto a smooth projective variety $Y$ with $\dim Y=n$. 
Let $L$ be a Cartier divisor on $X$ and let $\Delta$ be an 
effective $\mathbb R$-divisor on $X$ such that 
$K_X+\Delta$ is $\mathbb R$-Cartier. 
Let $D$ be an $\mathbb R$-divisor on $Y$. 
Let $k$ be a positive integer with $k\geq 2$. 
Assume the following conditions: 
\begin{itemize}
\item[(i)] $(X, \Delta)$ is log canonical 
over a nonempty Zariski open set of $Y$, and 
\item[(ii)] $L+f^*D-k(K_{X/Y}+\Delta)$ is nef and $f$-semi-ample. 
\end{itemize} 
Then the following properties hold. 
\begin{itemize}
\item[(1)] 
There exists a Cartier divisor $G$ on $Y$, 
which depends only on $f\colon X\to Y$, such that 
$$
\mathcal O_Y(G+\lceil D\rceil)\otimes f_*\mathcal O_X(L)
$$ is generically generated by global sections. 
\item[(2)] 
Let $H$ be a big Cartier divisor on $Y$ such that 
$|H|$ is free. 
Then $$\mathcal O_Y(K_Y+\lceil D\rceil+(n+1)H)\otimes 
\left(f_*\mathcal O_X(L)\right)^{**}$$ 
is generically generated by global sections. 
\item[(3)] 
Let $H^\dag$ be a nef and big Cartier divisor on $Y$ such that 
$|H^\dag|$ is not necessarily free. 
Then the sheaf 
$$\mathcal O_Y(K_Y+\lceil D\rceil+lH^\dag)\otimes 
\left(f_*\mathcal O_X(L)\right)^{**}$$ 
is generically generated by global sections 
for $l\geq n^2+2$. 
\end{itemize}
\end{cor}

\begin{proof} 
In Step \ref{p-step9.5-1}, we will treat (2) and (3), which are 
direct consequences of Theorem \ref{p-thm9.3}. 
In Step \ref{p-step9.5-2}, we will prove (1), which is much more 
difficult than (2) and (3). 

\setcounter{step}{0}
\begin{step}\label{p-step9.5-1} 
By Kodaira's lemma, we have $H\sim _{\mathbb Q}A+B$ such that 
$A$ is an ample $\mathbb Q$-divisor and $B$ is an effective 
$\mathbb Q$-divisor on $Y$. 
Let us consider 
$$
L+f^*H+f^*D-k(K_{X/Y}+\Delta)-f^*\left(\frac{1}{2}A+B\right). 
$$ 
Note that it is semi-ample by Lemma \ref{p-lem3.5}. 
We also note that $\frac{1}{2}A+B$ is big. 
Therefore, by Theorem \ref{p-thm9.3}, 
$$
\mathcal O_Y(K_Y+\lceil D\rceil +H)\otimes 
f_*\mathcal O_X(L) 
$$ 
is a mixed-$\widehat \omega$-big-sheaf on $Y$. 
Thus, by Lemma \ref{p-lem7.8}, 
$$
\mathcal O_Y(K_Y+\lceil D\rceil +(n+1)H)\otimes 
\left(f_*\mathcal O_X(L)\right)^{**}
$$ 
is generically generated by global sections. 
By the same argument, 
we see that 
$$
\mathcal O_Y(K_Y+\lceil D\rceil +H^\dag)\otimes 
f_*\mathcal O_X(L)
$$ 
is a mixed-$\widehat\omega$-big-sheaf on $Y$. 
Thus, by Lemma \ref{p-lem7.9}, 
the sheaf 
$$
\mathcal O_Y(K_Y+\lceil D\rceil +lH^\dag)\otimes 
\left(f_*\mathcal O_X(L)\right)^{**}
$$ 
is generically generated by global sections for $l\geq n^2+2$. 
\end{step}
\begin{step}\label{p-step9.5-2}
By taking a resolution of singularities as in Step 
\ref{p-step9.1.1} in the proof of Lemma \ref{p-lem9.1}, we may 
assume that $X$ is smooth. 
By replacing $L$ and $\Delta$ with 
$L+f^*\lceil D\rceil$ and $\Delta+\frac{1}{k}f^*\{-D\}$ 
respectively, we may assume that 
$D=0$. 
By the flattening theorem, there is a birational 
morphism $\tau\colon Y'\to Y$ from a smooth projective 
variety $Y'$ such that the main component 
of $X\times _Y Y'$ is flat over $Y'$. 
Let $X'$ be a resolution of the main component 
of $X\times _Y Y'$. 
Then we get the following commutative diagram. 
$$
\xymatrix{
X' \ar[d]_-{f'}\ar[r]^-\rho& X \ar[d]^-f\\
Y' \ar[r]_-\tau&Y
}
$$ 
By construction, any $f'$-exceptional 
divisor is $\rho$-exceptional. 
We put $K_{X'}+B=\rho^*(K_X+\Delta)$. 
We may assume that $\Supp B$ is a simple 
normal crossing divisor on $X'$. We write 
\begin{equation}\label{p-9.4}
K_{Y'}=\tau^*K_Y+R
\end{equation}
where $R$ is an effective $\tau$-exceptional 
divisor on $Y'$. We put 
$$
L':=\rho^*L+k\lceil -(B^{<0})\rceil -kf'^*R
$$ 
and $\Delta'=B+\lceil -(B^{<0})\rceil$. 
Note that $\lceil -(B^{<0})\rceil$ is effective and 
$\rho$-exceptional. 
Then we have 
$$
L'-k(K_{X'/Y'}+\Delta')=\rho^*(L-k(K_{X/Y}+\Delta)). 
$$ 
We take an effective $f'$-exceptional divisor 
$E$ on $X'$ such that 
$$
\left(f'_*\mathcal O_{X'}(L')\right)^{**} 
\simeq f'_*\mathcal O_{X'}(L'+E). 
$$ 
Note that $E$ is $\rho$-exceptional 
and that there is a generically isomorphic 
injection 
\begin{equation*}
\tau_*f'_*\mathcal O_{X'}(L'+E)
=f_*\rho_*\mathcal O_{X'}(L'+E)\subset 
f_*\mathcal O_X(L). 
\end{equation*} 
Therefore, we have a generically isomorphic 
injection 
\begin{equation}\label{p-9.5}
\tau_*\left((f'_*\mathcal O_{X'}(L'))^{**}\right)
\subset f_*\mathcal O_X(L). 
\end{equation} 
By Kodaira's lemma, we have $\tau^*H\sim _{\mathbb Q} A+B$ such that 
$A$ is an ample $\mathbb Q$-divisor and $B$ is an effective 
$\mathbb Q$-divisor. 
Note that 
\begin{equation*}
\begin{split}
&L'+E+f'^*\tau^*H-k\left(K_{X'/Y'}+\Delta'+\frac{1}{k}E+
\frac{1}{k}f'^*B\right)-\frac{1}{2}f'^*A\\
&=L'-k(K_{X'/Y'}+\Delta')+\frac{1}{2}f'^*A
\end{split}
\end{equation*}
is semi-ample by Lemma \ref{p-lem3.5}. 
Therefore, by Theorem \ref{p-thm9.3}, 
$$
\mathcal O_{Y'} (K_{Y'}+\tau^*H)\otimes 
f'_*\mathcal O_{X'}(L'+E)
$$ 
is a mixed-$\widehat\omega$-big-sheaf on $Y'$. 
Thus, by Lemma \ref{p-lem7.8}, 
$$
\mathcal O_{Y'}(K_{Y'}+(n+1)\tau^*H)
\otimes (f'_*\mathcal O_{X'}(L'))^{**}
$$ 
is generically generated by global sections. 
If we take a Cartier divisor $G$ on $Y$ such that 
$|\tau^*G-(K_{Y'}+(n+1)\tau^*H)|\ne \emptyset$, 
then $$\mathcal O_{Y'}(\tau^*G)\otimes 
\left(f'_*\mathcal O_{X'}(L')\right)^{**}$$ is generically generated by 
global sections. 
By \eqref{p-9.5}, we obtain that 
so is $\mathcal O_Y(G)\otimes f_*\mathcal O_X(L)$. 
\end{step}
We complete the proof of Corollary \ref{p-cor9.5}. 
\end{proof}

We note that \cite[Chapter V, 3.37.~Corollary]{nakayama} needs 
the assumption that 
$(X, \Delta)$ is klt over a nonempty Zariski open 
set of $Y$. On the other hand, 
Corollary \ref{p-cor9.5} can be applied 
to log canonical pairs. 
This is the main difference between 
\cite[Chapter V, 3.37.~Corollary]{nakayama} 
and Corollary \ref{p-cor9.5}. 

\medskip 

Sho Ejiri pointed out the following example, which was 
constructed by Hiroshi Sato. For some related example, 
see Example \ref{p-ex10.1} below. 

\begin{ex}[{\cite[Example 4.6]{fujino-gongyo}}]\label{p-ex9.6} 
There exists a flat toric morphism $f\colon X\to Y$ from a smooth 
projective toric threefold $X$ onto $Y=\mathbb P_{\mathbb P^1}
(\mathcal O_{\mathbb P^1}\oplus \mathcal O_{\mathbb P^1}(3))$. 
Let $\Delta$ be the union of all torus invariant 
divisors on $X$. 
Then it is well known that $(X, \Delta)$ is log canonical 
with $K_X+\Delta\sim 0$. 
In this case, 
$$
f_*\mathcal O_X(k(K_{X/Y}+\Delta))\otimes \mathcal O_Y(K_Y)
\simeq \mathcal O_Y(-(k-1)K_Y)
$$ 
holds for every integer $k$. 
Note that $-K_Y$ is not nef by 
$Y=\mathbb P_{\mathbb P^1}
(\mathcal O_{\mathbb P^1}\oplus \mathcal O_{\mathbb P^1}(3))$. 
Hence there are no ample Cartier divisors $A$ on $Y$ such that 
$$
f_*\mathcal O_X(k(K_{X/Y}+\Delta))\otimes \mathcal O_Y(K_Y+A)
\simeq \mathcal O_Y(-(k-1)K_Y+A)
$$ is generated by 
global sections for every positive integer $k\geq 2$. 
We note that if $\mathcal O_Y(-(k-1)K_Y+A)$ 
were generated by global sections for every integer $k\geq 2$ 
then $-K_Y$ would be nef.  
This is a contradiction. 
Therefore, we can not replace generic generations with 
global generations in Corollary \ref{p-cor9.5}. 
\end{ex}

\section{Proof of Theorems \ref{q-thm1.7}, 
\ref{q-thm1.8}, and \ref{q-thm1.9}}\label{p-sec10}

In this section, we prove Theorems \ref{q-thm1.7}, \ref{q-thm1.8}, 
and \ref{q-thm1.9} in Section \ref{q-sec1}. 

\medskip 

Let us first prove Theorem \ref{q-thm1.8}. 

\begin{proof}[Proof of Theorem \ref{q-thm1.8}]
We divide the proof into small steps. 
\setcounter{step}{0}
\begin{step}\label{p-step10-1}
By taking a suitable resolution of singularities of 
$X$, we may assume that 
$X$ is a smooth projective variety and $\Supp \Delta$ is 
a simple normal crossing divisor on $X$ (see Step \ref{p-step9.1.1} 
in the proof of Lemma \ref{p-lem9.1}). 
We may further assume that every log 
canonical center of $(X, \Delta_{\mathrm{hor}})$ is dominant onto $Y$. 
\end{step}

\begin{step}\label{p-step10-2}
In this step, we will prove the generic generation of 
$f_*\mathcal O_X(L)\otimes 
\mathcal O_Y(K_Y+lH)$ when $k=1$. 

\medskip

By replacing $L$ and $\Delta$ with $L-\lfloor 
\Delta_{\mathrm{ver}}\rfloor$ and $\Delta-
\lfloor 
\Delta_{\mathrm{ver}}\rfloor$ respectively, 
we may further assume that $(X, \Delta)$ is dlt and that every log 
canonical center of $(X, \Delta)$ is dominant onto $Y$. 
By the arguments in Step \ref{p-step7.7.2} in the 
Proof of Lemma \ref{p-lem7.7}, 
we see that $f_*\mathcal O_X(L)\otimes 
\mathcal O_Y(K_Y+lH)$ is generically generated by global sections. 
\end{step}

\begin{step}
In this step, 
we will see that 
$\left(f_*\mathcal O_X(L)\right)^{**}\otimes \mathcal O_Y(K_Y+lH)$ 
is generically generated by global sections when $k\geq 2$. 

\medskip 

This follows directly from Corollary \ref{p-cor9.5}. 
More precisely, we put $D=0$ and apply Corollary 
\ref{p-cor9.5} (2). 
\end{step}

\begin{step}
In this final step, we treat the case when $s\geq 2$. 
We take the $s$-fold fiber product 
$$
X^s:=\underbrace{X\times _Y X\times _Y \cdots \times _Y X}_{s}
$$ 
of $X$ over $Y$. 
Let $f^s\colon  X^s\to Y$ be the induced morphism. 
Let $\rho\colon X^{(s)}\to X^s$ be 
a resolution of singularities 
of the dominant components of $X^s$ such that 
$\rho$ is an isomorphism 
over a nonempty Zariski open set of $Y$. 
We put $f^{(s)}=f^s\circ \rho\colon  X^{(s)}\to Y$. 
We note that $X^{(s)}$ may be reducible, that is, 
a disjoint union of some smooth projective varieties. 
We can take a Zariski open set $U$ of $Y$ such that 
$\mathrm{codim}_Y(Y\setminus U)\geq 2$, 
$f_*\mathcal O_X(L)$ is locally free on $U$, 
and $f$ is flat over $U$. 
By applying Lemma \ref{p-lem3.7} to $f^{-1}(U)\to U$, we can construct 
a Cartier divisor $L^{(s)}$ on $X^{(s)}$ and an effective 
$\mathbb R$-divisor $\Delta^{(s)}$ on $X^{(s)}$ such that 
$$
L^{(s)}\sim_{\mathbb R} k(K_{X^{(s)}/Y}+\Delta^{(s)}), 
$$ 
$(X^{(s)}, \Delta^{(s)})$ is log canonical over a nonempty Zariski 
open set of $Y$, and there exists a generically isomorphic 
injection 
$$
\left(f^{(s)}_*\mathcal O_{X^{(s)}}(L^{(s)})\right)^{**}\subset 
\left(\bigotimes ^s f_*\mathcal O_X(L)\right)^{**}. 
$$ 
By Theorem \ref{p-thm9.3}, 
$$
\mathcal O_Y(K_Y+H)\otimes f^{(s)}_*\mathcal O_{X^{(s)}}(L^{(s)})
$$ 
is a finite direct sum of mixed-$\widehat\omega$-big-sheaves 
when $k\geq 2$. 
Note that $X^{(s)}$ may be reducible. 
Therefore, 
$$
\mathcal O_Y(K_Y+H)\otimes \left(\bigotimes ^s 
f_*\mathcal O_X(L)\right)^{**}
$$ 
is also a finite direct sum of mixed-$\widehat\omega$-big-sheaves. 
Thus, by Lemma \ref{p-lem7.8}, 
$$
\mathcal O_Y(K_Y+lH)\otimes \left(\bigotimes ^s 
f_*\mathcal O_X(L)\right)^{**}
$$
is generically generated by global sections 
for $l\geq n+1$ when $k\geq 2$ (see also Corollary \ref{p-cor9.5} (2)). 

If $k=1$, then we can check that 
$\mathcal O_Y(K_Y+lH)\otimes 
f^{(s)}_*\mathcal O_{X^{(s)}}(L^{(s)})$ is generically 
generated by global sections for $l\geq n+1$ by 
the arguments in Steps \ref{p-step10-1} and \ref{p-step10-2}. 
Therefore, 
$$
\mathcal O_Y(K_Y+lH)\otimes 
\left(\bigotimes ^s f_*\mathcal O_X(L)\right)^{**}
$$ 
is generically  generated by global sections for $l\geq n+1$ when 
$k=1$. 
\end{step} 
Hence we have obtained the desired statements. 
\end{proof}

Next we prove Theorem \ref{q-thm1.9}. 

\begin{proof}[Sketch of Proof of Theorem \ref{q-thm1.9}]
It is not difficult to modify the proof of 
Theorem \ref{q-thm1.8}. 
\setcounter{step}{0}
\begin{step} 
In this step, we will treat the case when $k=1$. 

\medskip 

As usual, by taking a suitable birational modification of $X$, 
we may assume that $X$ is smooth and $\Supp \Delta$ is 
a simple normal crossing divisor on $X$. 
By replacing $L$ and $\Delta$ with $L-\lfloor \Delta^{>1}\rfloor$ and 
$\Delta-\lfloor \Delta^{>1}\rfloor$ respectively, 
we may assume that $\Delta$ is a boundary $\mathbb R$-divisor 
on $X$. Note that 
$\Delta^{>1}$ is $f$-vertical. 
By perturbing 
the coefficients of $\Delta$, 
we may further assume that $\Delta$ is 
a $\mathbb Q$-divisor 
with $L\sim _{\mathbb Q} K_{X/Y}+\Delta$. 
By Lemma \ref{p-lem5.10}, 
$$
\mathcal O_X(K_X+\lfloor \Delta\rfloor 
+\lceil L-K_{X/Y}-\Delta\rceil)\simeq 
\mathcal O_X(L)\otimes f^*\mathcal O_Y(K_Y)
$$ 
is a mixed-$\omega$-sheaf on $X$. 
Therefore, $f_*\mathcal O_X(L)\otimes 
\mathcal O_Y(K_Y)$ is a mixed-$\omega$-sheaf on $Y$. 
Thus, by Lemma \ref{p-lem7.9}, 
$f_*\mathcal O_X(L)\otimes \mathcal O_Y(K_Y+lH^\dag)$ is 
generically generated by global sections for $l\geq n^2+1$. 
Similarly, 
we may assume that the sheaf 
$f^{(s)}_*\mathcal O_{X^{(s)}}(L^{(s)})\otimes \mathcal O_Y(K_Y)$ 
in the proof of Theorem \ref{q-thm1.8} is a finite direct sum of 
mixed-$\omega$-sheaves on $Y$ when $k=1$. 
Therefore, 
$$
\left(\bigotimes ^sf_*\mathcal O_X(L)\right)^{**}
\otimes \mathcal O_Y(K_Y+lH^\dag)
$$ 
is generically generated by global sections for $l\geq n^2+1$. 
\end{step}
\begin{step}
In this step, we will treat the case when $k\geq 2$. 

\medskip 

If we use Lemma \ref{p-lem7.9} instead of Lemma \ref{p-lem7.8}, 
then the proof of Theorem \ref{q-thm1.8} implies that 
$$
\mathcal O_Y(K_Y+lH^\dag)\otimes 
\left(\bigotimes ^sf_*\mathcal O_X(L)\right)^{**}
$$ 
is generically generated by global sections 
for $l\geq n^2+2$ (see also Corollary \ref{p-cor9.5} (3)). 
\end{step}
Thus we get the desired statements. 
\end{proof}

Finally, we prove Theorem \ref{q-thm1.7}. 

\begin{proof}[Proof of Theorem \ref{q-thm1.7}]
We put $L=kK_{X/Y}$. Then this theorem directly follows 
from Theorems \ref{q-thm1.8} and 
\ref{q-thm1.9}. 
\end{proof}

In \cite[Section 8]{fujino-fujisawa}, we constructed 
the following example, 
which shows that we can not replace the 
generic generation with 
the global generation in Conjecture \ref{q-conj1.5}.  

\begin{ex}\label{p-ex10.1}
There exists a surjective morphism $f\colon X\to Y$ between smooth 
projective varieties with the following properties. 
\begin{itemize}
\item[(i)] $Y$ is a Kummer surface. In particular, 
$\omega_Y\simeq \mathcal O_Y$ holds. 
\item[(ii)] $C_i$ is a $(-2)$-curve on $Y$ for $1\leq i\leq 16$. 
\item[(iii)] $f$ is smooth over $U=Y\setminus \sum _{i=1}^{16}C_i$. 
\item[(iv)] $L_{X/Y}$ is a Weil divisor on $Y$ which is numerically 
equivalent to $\frac{1}{2}\sum _{i=1}^{16}C_i$. 
\item[(v)] We have 
$$
f_*\omega^{\otimes k}_{X/Y} \simeq 
\begin{cases}
\mathcal O_Y(\sum _{i=1}^{16} lC_i) & k=2l, \\ 
\mathcal O_Y(L_{X/Y}+\sum _{i=1}^{16} lC_i) &k=2l+1. 
\end{cases}
$$
\end{itemize}
Let $H$ be an ample Cartier divisor on $Y$. 
By Reider's theorem (see, for example, \cite[Chapter 
IV, (11.4) Theorem]{bhpv}), 
$|3H|$ is free. 
We note that $L_{X/Y}+3H$ is ample by Nakai's ampleness 
criterion. 
Therefore, by Reider's theorem again, $|L_{X/Y}+3H|$ is 
free. 
This means that 
$$
f_*\omega^{\otimes k}_{X/Y}\otimes 
\omega_Y\otimes \mathcal O_Y(3H)
$$ is generated by 
global sections on $U$. 
On the other hand, we have 
$$
\left(f_*\omega^{\otimes k}_{X/Y}\otimes 
\omega_Y\otimes \mathcal O_Y(3H)\right)\cdot C_i=-k+3H\cdot C_i. 
$$ 
Therefore, if $k>3H\cdot C_{i_0}$ holds for 
some $1\leq i_0\leq 16$, 
then 
$$
f_*\omega^{\otimes k}_{X/Y}\otimes 
\omega_Y\otimes \mathcal O_Y(3H)
$$ 
is not generated by global sections. 
\end{ex}

We close this section with an easy remark. 

\begin{rem}\label{p-rem10.2}
Let $Y$ be a smooth projective variety and let $H$ be an ample 
Cartier divisor on $Y$. Let $m$ be any positive integer. 
Then we can construct a finite cover $f\colon X\to Y$ from 
a smooth projective variety $X$ such that 
$\mathcal O_Y(-mH)$ is a direct summand of $f_*\mathcal O_X$. 
Therefore, we need the condition $k\geq 1$ 
in Theorems \ref{q-thm1.7}, \ref{q-thm1.8}, and 
\ref{q-thm1.9}. 
\end{rem}

\section{Some other applications}\label{p-sec11}

In this section, we treat Nakayama's inequality on $\kappa _\sigma$ and 
a slight generalization of 
the twisted weak positivity theorem. 
Theorem \ref{p-thm11.3} and a special case of 
Theorem \ref{p-thm11.7} have already played a crucial 
role in the theory of minimal models. 

\medskip 

Let us first recall the definition of $\kappa _\sigma$ for the 
reader's convenience. 

\begin{defn}[{Nakayama's numerical dimension, 
see \cite[Chapter V.2.5.~Definition]{nakayama}}]\label{p-def11.1} 
Let $D$ be a pseudo-effective $\mathbb R$-Cartier divisor 
on a normal projective variety $X$ and let $A$ be a Cartier divisor on $X$. 
If $H^0(X, \mathcal O_X(\lfloor mD\rfloor +A))\ne 0$ for infinitely 
many positive integers $m$, then we set 
$$
\sigma (D; A)
=\max \left\{ k\in \mathbb Z_{\geq 0}\, \left|\, {\underset{m\to \infty}{\limsup}}
\frac{\dim H^0(X, \mathcal O_X(\lfloor mD\rfloor +A))}{m^k}>0 \right.\right\}. 
$$
If $H^0(X, \mathcal O_X(\lfloor mD\rfloor +A))\ne 0$ only for 
finitely many $m\in \mathbb Z_{\geq 0}$, then 
we set $\sigma (D; A)=-\infty$. 
We define {\em{Nakayama's numerical dimension $\kappa _{\sigma}$}} 
by
$$
\kappa _\sigma(X, D)=
\max \{\sigma(D; A)\, |\, A \ {\text{is a  Cartier divisor on}} \ X\}. 
$$
It is well known that $\kappa_{\sigma}(X, D)\geq 0$ 
(see, for example, \cite[Chapter V.~2.7.~Proposition]{nakayama}). 
If $D$ is not pseudo-effective, then we put $\kappa_\sigma
(X, D)=-\infty$. By this convention, we can define 
$\kappa_\sigma(X, D)$ 
for every $\mathbb R$-Cartier divisor $D$ on $X$. 
It is obvious that 
$$
\kappa_\sigma(X, D)\geq \kappa(X, D) 
$$
always 
holds for every $\mathbb R$-Cartier divisor $D$ on $X$ by definition, 
where $\kappa (X, D)$ denotes 
the {\em{Iitaka dimension}} of $D$. 
\end{defn}

For the details of $\kappa _\sigma(X, D)$ and $\kappa (X, D)$, 
we recommend the reader to see \cite{nakayama}. 
The following remark is easy but very useful. 

\begin{rem}[{\cite[Chapter V, 2.6.~Remark (6)]{nakayama}}]
\label{p-rem11.2}
Let $X$ be a smooth projective variety and let $D$ be an $\mathbb 
R$-divisor on $X$. 
We put 
$$
\sigma (D; A)'=\max \left\{ k\in \mathbb Z_{\geq 0}\cup 
\{-\infty\}\, \left|\, {\underset{m\to \infty}{\limsup}}
\frac{\dim H^0(X, \mathcal O_X(\lceil mD\rceil +A))}{m^k}>0 \right.\right\}, 
$$ 
where $A$ is a divisor on $X$. 
Then we have the following equality 
$$
\kappa _\sigma(X, D)=\max\{ \sigma(D; A)' \, |\, 
\text{$A$ is a divisor}\}. 
$$
We will use this characterization of $\kappa_\sigma$ in the proof of 
Theorem \ref{p-thm11.3} below. 

We note the following 
easy but important fact 
that $\kappa _\sigma(X, lD)=\kappa_\sigma(X, D)$ holds for every positive 
integer $l$ (see \cite[Remark 2.2]{fujino-corri}), which 
will be useful in the proof of Theorem \ref{p-thm11.3} below. 
\end{rem}

The inequalities in Theorem \ref{p-thm11.3} are indispensable 
in the theory of minimal models (see 
Remarks \ref{p-rem11.4} and \ref{p-rem11.5}). 

\begin{thm}[{\cite[Chapter V, 4.1.~Theorem (1)]{nakayama} and 
\cite[Section 3]{fujino-corri}}]\label{p-thm11.3} 
Let $f\colon X\to Y$ be a surjective morphism from a normal 
projective variety $X$ onto a smooth projective variety $Y$ 
with connected fibers. 
Let $\Delta$ be an effective $\mathbb R$-divisor 
on $X$ such that 
$K_X+\Delta$ is $\mathbb R$-Cartier and that 
$(X, \Delta)$ is log canonical over a nonempty Zariski open set 
of $Y$. 
Let $D$ be an $\mathbb R$-Cartier $\mathbb R$-divisor 
on $X$ such that $D-(K_{X/Y}+\Delta)$ is nef. 
Then, for any $\mathbb R$-divisor $Q$ on $Y$, 
we have 
$$
\kappa _\sigma(X, D+f^*Q)\geq \kappa _\sigma(F, D|_F)
+\kappa (Y, Q)
$$ 
and 
$$
\kappa _\sigma(X, D+f^*Q)\geq \kappa (F, D|_F)
+\kappa _\sigma(Y, Q)
$$ 
where $F$ is a sufficiently general fiber of $f\colon X\to Y$. 
\end{thm}

Before we prove Theorem \ref{p-thm11.3}, 
we give two important remarks. 

\begin{rem}\label{p-rem11.4}
We think that one of the most important results of 
Nakayama's theory of $\omega$-sheaves is the 
inequality on $\kappa_\sigma$ in 
\cite[Chapter V, 4.1.~Theorem (1)]{nakayama}. 
However, as we explained in \cite[Remark 3.8]{fujino-subadd} and 
\cite[Section 3]{fujino-corri}, 
the proof of 
\cite[Chapter V, 4.1.~Theorem (1)]{nakayama} is 
incomplete. For the details, see, for example, \cite[Section 1]{fujino-corri}. 
So, in Theorem \ref{p-thm11.3}, 
we claim two weaker inequalities than Nakayama's original 
one (see \cite[(3.3) and (3.4)]{fujino-corri}). 
The first inequality in Theorem \ref{p-thm11.3} 
is still sufficiently powerful for some 
geometric applications (see \cite[Section 3]{fujino-corri}). 
\end{rem}

\begin{rem}[{see \cite[Section 3]{fujino-corri}}]
\label{p-rem11.5} 
The troubles 
in the proof of \cite[Remark 2.6]{dhp} and 
\cite[Theorem 4.3]{gongyo-lehmann} 
caused by the incompleteness of 
\cite[Chapter V, 4.1.~Theorem (1)]{nakayama} can be 
corrected by using the first inequality in Theorem \ref{p-thm11.3}. 
For the details, we recommend the reader to see 
\cite[Lemma 2.11]{hashizume-hu}. 
\end{rem}

Let us prove Theorem \ref{p-thm11.3}. 

\begin{proof}[Proof of Theorem \ref{p-thm11.3}] 
If $Q$ is not pseudo-effective, then the desired inequalities 
are obviously true. So we may assume that $Q$ is pseudo-effective. 
Similarly, we may further assume that $D|_F$ is pseudo-effective. 
As usual (see Step \ref{p-step9.1.1} in the proof of 
Lemma \ref{p-lem9.1}), 
we may assume that $X$ is smooth and 
$\Supp \Delta$ is a simple normal crossing divisor on $X$ 
by the basic properties of $\kappa _\sigma$ and 
$\kappa$. 
We take a sufficiently ample Cartier divisor $A$ on $X$ such that 
$A+\{-mD\}$ is ample for every integer $m$. 
Then 
$$
\lceil mD\rceil +A-m(K_{X/Y}+\Delta)=m(D-(K_{X/Y}+\Delta))+A+\{-mD\}
$$ 
is ample for every positive integer $m$. 
Then we can take an ample Cartier divisor 
$H$ on $Y$ such that 
$\mathcal O_Y(H)\otimes f_*\mathcal O_X(\lceil mD\rceil +A)$ 
is generically generated by global sections for every positive 
integer $m$ by Corollary \ref{p-cor9.5} (1). 
Thus there exists a generically isomorphic injection 
$$
\mathcal O_Y^{\oplus r(mD; A)}\hookrightarrow 
\mathcal O_Y(H)\otimes f_*\mathcal O_X(\lceil mD\rceil +A), 
$$
where $r(mD; A):=\rank f_*\mathcal O_X(\lceil mD\rceil +A)$. 
This induces the following injection 
$$
\mathcal O_Y(\lfloor mQ\rfloor +H)^{\oplus 
r(mD; A)}\hookrightarrow 
\mathcal O_Y(\lfloor mQ\rfloor +2H)\otimes f_*\mathcal O_X(\lceil mD\rceil 
+A). 
$$ 
Therefore, we have 
\begin{equation}\label{p-11.1}
\begin{split}
&\dim _{\mathbb C} H^0(X, \mathcal O_X(\lceil m(D+f^*Q)\rceil 
+A+2f^*H)
\\&\geq \dim _{\mathbb C} H^0(X, \mathcal O_X(\lceil mD\rceil+f^*(\lfloor 
mQ\rfloor)+A+2f^*H)) 
\\ & \geq 
r(mD; A)\cdot \dim _{\mathbb C}
H^0(Y, \mathcal O_Y(\lfloor mQ\rfloor +H))
\end{split}
\end{equation} 
for every positive integer $m$. 
We can take a positive integer $m_0$ and a 
positive real number $C_0$ such that 
\begin{equation}\label{p-11.2}
C_0m^{\kappa (F, D|_F)}\leq r(mm_0D; A)
\end{equation}
for every large positive integer $m$ (see, 
for example, \cite[Chapter II, 3.7.~Theorem]
{nakayama}). 
Thus we have 
\begin{equation}\label{p-11.3}
\begin{split}
&\dim H^0(X, \mathcal O_X(\lceil mm_0(D+f^*Q)\rceil+A+2f^*H))\\
&\geq C_0m^{\kappa(F, D|_F)}
\cdot \dim H^0(Y, \mathcal O_Y(\lfloor mm_0Q\rfloor +H)) 
\end{split}
\end{equation} 
for every large positive integer $m$ by \eqref{p-11.1} and 
\eqref{p-11.2}. 
We may assume that $H$ is sufficiently ample. 
Then we get 
\begin{equation}\label{p-11.4}
\limsup_{m\to \infty} \frac{\dim H^0(X, 
\mathcal O_X(\lceil mm_0(D+f^*Q)\rceil+A+2f^*H))}
{m^{\kappa (F, D|_F)+\kappa _{\sigma}(Y, Q)}}>0
\end{equation}
by \eqref{p-11.3} and the definition of $\kappa _\sigma(Y, Q)$. 
This means that the following inequality 
\begin{equation}\label{p-11.5}
\kappa _\sigma(X, D+f^*Q)\geq \kappa (F, D|_F)+
\kappa _\sigma(Y, Q)
\end{equation} 
holds. 

Similarly, we can take a positive integer $m_1$ and 
a positive real number $C_1$ such that 
\begin{equation}\label{p-11.6}
\begin{split}
C_1m^{\kappa (Y, Q)}&\leq 
\dim H^0(Y, \mathcal O_Y(\lfloor mm_1Q\rfloor))\\ 
&\leq \dim H^0(Y, \mathcal O_Y(\lfloor mm_1Q\rfloor+H))
\end{split}
\end{equation}
for every large positive integer $m$ (see, 
for example, 
\cite[Chapter II, 3.7.~Theorem]
{nakayama}) if 
$H$ is a sufficiently ample Cartier divisor. 
Then, by \eqref{p-11.1} and \eqref{p-11.6}, we have 
\begin{equation}\label{p-11.7}
\begin{split}
&\dim H^0(X, \mathcal O_X(\lceil mm_1 (D+f^*Q)\rceil 
+A+2f^*H))\\
&\geq C_1m^{\kappa(Y, Q)}
\cdot r(mm_1D; A)
\end{split}
\end{equation}
for every large positive integer $m$. 
Therefore, we get 
\begin{equation}\label{p-11.8}
\limsup_{m\to \infty} \frac{\dim H^0(X, 
\mathcal O_X(\lceil mm_1 (D+f^*Q)\rceil +A+2f^*H))}
{m^{\kappa _\sigma(F, D|_F)+\kappa(Y, Q)}}>0 
\end{equation}
when $A$ is sufficiently ample. 
Note that 
\begin{equation}\label{p-11.9}
\sigma (m_1D|_F; A|_F)'=\max \left\{
k\in \mathbb Z_{\geq 0} \cup \{-\infty\}\, 
\left|\, \underset{m\to \infty}{\limsup}\frac{r(mm_1D; A)}{m^k}>0\right.\right\} 
\end{equation}
for a sufficiently general fiber $F$ of $f\colon X\to Y$ and 
that 
\begin{equation}\label{p-11.10}
\begin{split}
\kappa_\sigma(F, D|_F)
&=\kappa_\sigma (F, m_1D|_F)
\\&=\max\{\sigma(m_1D|_F; A|_F)'\, |\, {\text{$A$ is very ample}}\}. 
\end{split}
\end{equation}
Hence we have the inequality 
\begin{equation}\label{p-11.11}
\kappa _\sigma(X, D+f^*Q)\geq \kappa_\sigma (F, D|_F)+
\kappa (Y, Q) 
\end{equation} 
by \eqref{p-11.8}. 
\end{proof}

It is highly desirable to solve the following conjecture. 
As we explained in \cite{fujino-corri}, 
Nakayama's original inequality on $\kappa_\sigma$ 
(see \cite[Chapter V, 4.1.~Theorem (1)]{nakayama}) follows 
from Conjecture \ref{p-conj11.6} and the argument in the proof of 
Theorem \ref{p-thm11.3}. 

\begin{conj}[{\cite[Conjecture 1.4]{fujino-corri}}]\label{p-conj11.6}
Let $X$ be a smooth projective variety and let $D$ be a 
pseudo-effective $\mathbb R$-divisor 
on $X$. 
Then there exist a positive integer $m_0$, a 
positive rational number $C$, and an ample Cartier divisor 
$A$ on $X$ such that 
\begin{equation*}
Cm^{\kappa_\sigma(X, D)}\leq \dim H^0(X, \mathcal O_X
(\lfloor mm_0D\rfloor +A))
\end{equation*}
holds for every large positive integer $m$. 
\end{conj}

Finally, we treat a slight 
generalization of the twisted weak positivity theorem. 

\begin{thm}[Twisted weak positivity theorem]\label{p-thm11.7}
Let $f\colon X\to Y$ be a surjective morphism 
from a normal projective variety $X$ onto a smooth 
projective variety $Y$. 
Let $\Delta$ be an effective $\mathbb R$-divisor 
on $X$ such that $K_X+\Delta$ is $\mathbb R$-Cartier and 
that $(X, \Delta)$ is log canonical over a nonempty 
Zariski open set of $Y$. 
Let $L$ be a Cartier divisor on $X$ with $L\sim _{\mathbb R} 
k(K_{X/Y}+\Delta)$ for some 
positive integer $k$. 
Then the sheaf 
$
f_*\mathcal O_X(L)
$ 
is weakly positive. 
\end{thm}

\begin{proof}
Let $\alpha$ be a positive integer and let $\mathcal H$ 
be an ample invertible sheaf on $Y$. 
By Theorem \ref{q-thm1.8} or Theorem \ref{q-thm1.9}, 
we can take a positive integer $\beta$ which depends only 
on $Y$ such that 
$$
\left(\bigotimes ^s f_*\mathcal O_X(L)\right)^{**} 
\otimes \mathcal H^{\otimes \beta}
$$ 
is generically generated by global sections for every positive 
integer $s$. 
This implies that 
$$
\widehat S^{\alpha\beta} (f_*\mathcal O_X(L))\otimes 
\mathcal H^{\otimes \beta}
$$
is generically generated by global sections. 
This means that 
$f_*\mathcal O_X(L)$ is weakly positive. 
\end{proof}

\section{On Iwai's theorem:~Theorem \ref{q-thm1.6}}\label{p-sec12}

This section is independent of the other sections. 
Here we explain the following result due to Masataka Iwai. 
The proof of Theorem \ref{p-thm12.1} is analytic and is completely 
different from the arguments in this paper.  

\begin{thm}[Masataka Iwai]\label{p-thm12.1} 
Let $f\colon X\to Y$ be a surjective morphism between 
smooth projective varieties with connected fibers 
and let $\mathcal L$ be an ample invertible sheaf on $Y$. 
Let $U$ be the largest Zariski open set of $Y$ such that 
$f$ is smooth over $U$. 
We put $\dim Y=n$. 
Then 
$$
f_*\omega^{\otimes a}_{X/Y}\otimes \omega_Y\otimes \mathcal L^{\otimes b}
$$ 
is generated by global sections on $U$ for all integers 
$a\geq 1$ and $b\geq \frac{n(n+1)}{2}+1$. 
\end{thm}

\begin{proof}[Sketch of Proof] 
Here we will only explain how to modify the proof of 
\cite[Theorem 1.4]{iwai}. 
As in \cite[Theorem 2.3]{iwai}, 
we take a smooth hermitian metric $h_{\mathcal L}$ on $\mathcal L$, 
a K\"ahler metric $\omega$ on $Y$, and a quasi-plurisubharmonic 
function $\varphi$ on $Y$. 
Let $h_a$ be the singular hermitian 
metric on $\omega^{\otimes a}_{X/Y}$ in 
\cite[Theorem 2.4]{iwai}. We put 
$$
\mathcal L^\dag=\omega^{\otimes (a-1)}_{X/Y} 
\otimes f^*\mathcal L^{\otimes (N+1+b^\dag)} 
\quad \text{and}\quad 
h_{\mathcal L^\dag}=h^{\frac{a-1}{a}}_a\!f^*h^{N+1+b^\dag}
_{\mathcal L}, 
$$ 
where $N=\frac{n(n+1)}{2}$ and $b^\dag=b-(N+1)\geq 0$. 
Then we consider the adjoint bundle 
$$
\omega_X\otimes \mathcal L^\dag
\simeq \omega^{\otimes a}_{X/Y} 
\otimes f^*(\omega_Y\otimes \mathcal L^{\otimes b}). 
$$ 
In this situation, the proof of \cite[Theorem 1.4]{iwai} 
implies that 
$$
H^0(Y, f_*(\omega_X\otimes \mathcal L^\dag))\otimes 
\mathcal O_Y\to f_*(\omega_X\otimes \mathcal L^\dag)
$$ 
is surjective on $U$, equivalently, 
$$
H^0(Y, f_*\omega^{\otimes a}_{X/Y}\otimes \omega_Y\otimes 
\mathcal L^{\otimes b})\otimes 
\mathcal O_Y\to f_*\omega^{\otimes a}_{X/Y}\otimes \omega_Y\otimes 
\mathcal L^{\otimes b}
$$ 
is surjective on $U$. This is what we wanted. 
\end{proof}


\begin{thebibliography}{BHPV}

\bibitem[AK]{abramovich-karu} 
D.~Abramovich, K.~Karu, 
Weak semistable reduction in characteristic $0$, 
Invent. Math. \textbf{139} (2000), no. 2, 241--273. 

\bibitem[BHPV]{bhpv} 
W.~P.~Barth, K.~Hulek, C.~A.~M.~Peters, A.~Van de Ven, 
{\em{Compact complex surfaces}}, 
Second edition. Ergebnisse der Mathematik und 
ihrer Grenzgebiete. 3. Folge. A Series of Modern 
Surveys in Mathematics [Results in 
Mathematics and Related Areas. 3rd 
Series. A Series of Modern Surveys 
in Mathematics], \textbf{4}. Springer-Verlag, 
Berlin, 2004. 

\bibitem[C]{conrad} 
B.~Conrad, {\em{Grothendieck duality and base change}}, 
Lecture Notes in Mathematics, \textbf{1750}. Springer-Verlag, 
Berlin, 2000.

\bibitem[DHP]{dhp} 
J.-P.~Demailly, C.~D.~Hacon, M.~P\u aun, 
Extension theorems, non-vanishing and the existence of 
good minimal models, 
Acta Math. \textbf{210} (2013), no. 2, 203--259.

\bibitem[De]{denf}
Y.~Deng, Applications of the Ohsawa--Takegoshi extension theorem to 
direct image problems, Int. Math. Res. Not. IMRN 2021, 
no. 23, 17611--17633.

\bibitem[Du]{dutta} 
Y.~Dutta, On the effective freeness of the direct images of 
pluricanonical bundles, 
Ann. Inst. Fourier (Grenoble) \textbf{70} (2020), no. 4, 1545--1561.

\bibitem[DuM]{dutta-murayama} 
Y.~Dutta, T.~Murayama, Effective generation and twisted weak positivity 
of direct images, 
Algebra Number Theory \textbf{13} (2019), no. 2, 425--454.

\bibitem[EKL]{ekl} 
L.~Ein, O.~K\"uchle, R.~Lazarsfeld, 
Local positivity of ample line bundles, 
J. Differential Geom. \textbf{42} (1995), no. 2, 193--219.

\bibitem[EV]{esnault-viehweg} 
H.~Esnault, E.~Viehweg, 
{\em{Lectures on vanishing theorems}}, 
DMV Seminar, {\textbf{20}}. Birkh\"auser Verlag, Basel, 1992.

\bibitem[Fn1]{fujino-higher}
O.~Fujino, 
Higher direct images of log canonical divisors, 
J. Differential Geom. \textbf{66} (2004), no. 3, 453--479.

\bibitem[Fn2]{fujino-fundamental} 
O.~Fujino, 
Fundamental theorems for the log minimal model program, 
Publ. Res. Inst. Math. Sci. \textbf{47} (2011), no. 3, 727--789.

\bibitem[Fn3]{fujino-quasi-alb} 
O.~Fujino, On quasi-Albanese maps, preprint (2014). 

\bibitem[Fn4]{fujino-direct}
O.~Fujino, 
Direct images of relative pluricanonical bundles, 
Algebr. Geom. \textbf{3} (2016), no. 1, 50--62.

\bibitem[Fn5]{fujino-direct-corri}
O.~Fujino, 
Corrigendum:~Direct images of 
relative pluricanonical bundles (Algebraic Geometry \textbf{3}, no. 1, (2016), 50--62), 
Algebr. Geom. \textbf{3} (2016), no. 2, 261--263.

\bibitem[Fn6]{fujino-foundations} 
O.~Fujino, {\em{Foundations of the minimal model program}}, 
MSJ Memoirs, \textbf{35}. Mathematical Society of Japan, Tokyo, 2017. 

\bibitem[Fn7]{fujino-subadd}
O.~Fujino, On subadditivity of the logarithmic Kodaira dimension, 
J. Math. Soc. Japan \textbf{69} (2017), no. 4, 1565--1581. 

\bibitem[Fn8]{fujino-zucker} 
O.~Fujino, On semipositivity, injectivity and vanishing theorems, 
{\em{Hodge theory and $L^2$-analysis}}, 
245--282, Adv. Lect. Math. (ALM), \textbf{39}, Int. Press, 
Somerville, MA, 2017.

\bibitem[Fn9]{fujino-weak} 
O.~Fujino, Notes on the weak positivity theorems, 
{\em{Algebraic varieties and automorphism groups}}, 73--118, 
Adv. Stud. Pure Math., \textbf{75}, Math. Soc. Japan, Tokyo, 2017. 

\bibitem[Fn10]{fujino-iitaka}
O.~Fujino, {\em{Iitaka conjecture:~An introduction}}, 
SpringerBriefs in Mathematics. Springer, Singapore, 2020. 

\bibitem[Fn11]{fujino-corri} 
O.~Fujino, Corrigendum to \lq\lq On 
subadditivity of the logarithmic Kodaira dimension\rq\rq, 
J. Math. Soc. Japan \textbf{72} (2020), no. 4, 1181--1187. 

\bibitem[Fn12]{fujino-nakayama} 
O.~Fujino, On Nakayama's theorem, 
J. Math. Sci. Univ. Tokyo \textbf{28} (2021), no. 4, 641--650. 

\bibitem[FF]{fujino-fujisawa}
O.~Fujino, T.~Fujisawa, 
Variations of mixed Hodge structure and semipositivity theorems, 
Publ. Res. Inst. Math. Sci. \textbf{50} (2014), no. 4, 589--661. 

\bibitem[FFS]{ffs} 
O.~Fujino, T.~Fujisawa, M.~Saito, 
Some remarks on the semipositivity theorems, 
Publ. Res. Inst. Math. Sci. \textbf{50} (2014), no. 1, 85--112.

\bibitem[FG]{fujino-gongyo} 
O.~Fujino, Y.~Gongyo, On images of weak Fano manifolds, 
Math. Z. \textbf{270} (2012), no. 1-2, 531--544.

\bibitem[Fs]{fujisawa} 
T.~Fujisawa, A remark on semipositivity theorems, 
preprint (2017). arXiv:1710.01008 [math.AG]

\bibitem[GL]{gongyo-lehmann} 
Y.~Gongyo, B.~Lehmann, 
Reduction maps and minimal model theory, 
Compos. Math. \textbf{149} (2013), no. 2, 295--308. 

\bibitem[H]{hartshorne} 
R.~Hartshorne, 
{\em{Residues and duality}}, 
Lecture notes of a seminar on the work 
of A.~Grothendieck, given at Harvard 1963/64. With an 
appendix by P.~Deligne. Lecture Notes in 
Mathematics, No. \textbf{20} Springer-Verlag, Berlin-New York, 1966. 

\bibitem[HH]{hashizume-hu} 
K.~Hashizume, Z.~Hu, 
On minimal model theory for log abundant lc pairs, 
J. Reine Angew. Math. \textbf{767} (2020), 109--159.

\bibitem[I]{iwai} 
M.~Iwai, On the global generation of direct images of pluri-adjoint line bundles, 
Math. Z. \textbf{294} (2020), no. 1-2, 201--208.

\bibitem[Ka]{kawamata-curve} 
Y.~Kawamata, 
Kodaira dimension of algebraic fiber spaces over curves, 
Invent. Math. \textbf{66} (1982), no. 1, 57--71.

\bibitem[Ko]{kollar} 
J.~Koll\'ar, Higher direct images of dualizing sheaves. I, 
Ann. of Math. (2) \textbf{123} (1986), no. 1, 11--42.

\bibitem[KM]{kollar-mori} 
J.~Koll\'ar, S.~Mori, 
{\em{Birational geometry of algebraic varieties}}. With 
the collaboration of C.~H.~Clemens and A.~Corti. Translated 
from the 1998 Japanese original. Cambridge 
Tracts in Mathematics, {\textbf{134}}. Cambridge University 
Press, Cambridge, 1998.

\bibitem[M]{mori} 
S.~Mori, Classification of higher-dimensional 
varieties, {\em{Algebraic geometry, Bowdoin, 1985 (Brunswick, Maine, 1985)}}, 
269--331, Proc. Sympos. Pure Math., \textbf{46}, Part 1, 
Amer. Math. Soc., Providence, RI, 1987.

\bibitem[N]{nakayama} 
N.~Nakayama, {\em{Zariski-decomposition and abundance}}, 
MSJ Memoirs, {\textbf{14}}. Mathematical Society of Japan, Tokyo, 2004. 

\bibitem[PS]{popa-schnell}
M.~Popa, C.~Schnell, 
On direct images of pluricanonical bundles, 
Algebra Number Theory \textbf{8} (2014), no. 9, 2273--2295. 
 
\bibitem[SY]{su-yang} 
X.~Su, X.~Yang, 
Global generation and very ampleness for adjoint 
linear series, 
Comm. Anal. Geom. \textbf{27} (2019), no. 7, 1639--1663.  

\bibitem[V1]{viehweg1}
E.~Viehweg, Weak positivity and the additivity of the Kodaira 
dimension for certain fibre spaces, 
{\em{Algebraic varieties and analytic varieties (Tokyo, 1981)}}, 329--353, 
Adv. Stud. Pure Math., 1, North-Holland, Amsterdam, 1983.

\bibitem[V2]{viehweg2} 
E.~Viehweg, Weak positivity and the additivity of the Kodaira 
dimension. II. The local Torelli map, {\em{Classification of algebraic 
and analytic manifolds (Katata, 1982)}}, 567--589, 
Progr. Math., \textbf{39}, Birkh\"auser Boston, Boston, MA, 1983.

\bibitem[V3]{viehweg3} 
E.~Viehweg, 
{\em{Quasi-projective moduli for polarized manifolds}}, 
Ergebnisse der Mathematik und ihrer 
Grenzgebiete (3) [Results in Mathematics and 
Related Areas (3)], \textbf{30}. Springer-Verlag, Berlin, 1995.

\end{thebibliography}
\end{document}